\documentclass[a4paper,reqno]{amsart}

\usepackage{amsmath}
\usepackage{amsthm}
\usepackage{amsfonts}
\usepackage{amssymb}
\usepackage{float}

\newcounter{rh}
\newcounter{lp}
\usepackage{graphicx}
\usepackage{subcaption}

\usepackage{pdfsync}

\usepackage{asymptote}
\usepackage{verbatim}

\usepackage{subcaption}

\usepackage{color}

\usepackage{hyperref}

\usepackage{graphicx}

\hypersetup{
    colorlinks=false,
    pdfborder={0 0 0},
    pdfstartview={XYZ null null 1.00}
}

\numberwithin{equation}{section}

  \def\<{\langle}
  \def\>{\rangle}

  \def\ve{\varepsilon}

  \def\R{\mathbb{R}}

  \def\C{\mathbb{C}}

  \def\RR{\mathcal{R}}

  \def\CC{\mathcal{C}}
  \def\SS{\mathcal{S}}

  \def\t{\widetilde}

\theoremstyle{plain}
  \newtheorem{theorem}{Theorem}[section]

  \newtheorem{proposition}[theorem]{Proposition}
  
  \newtheorem{lemma}[theorem]{Lemma}

\theoremstyle{definition}

\begin{document}

\title[On global solutions of defocusing mKdV equation...]{On global solutions of defocusing mKdV equation with specific initial data of critical regularity}

\author{Kamil Dunst}
\address{\noindent Faculty of Mathematics and Computer Science \newline Nicolaus Copernicus University \newline Chopina 12/18, 87-100 Toru\'n, Poland}
\email{globsztajn@mat.umk.pl}

\author{Piotr Kokocki}
\address{\noindent Faculty of Mathematics and Computer Science \newline Nicolaus Copernicus University \newline Chopina 12/18, 87-100 Toru\'n, Poland}
\email{pkokocki@mat.umk.pl}
\thanks{The researches supported by the MNiSW Iuventus Plus Grant no. 0338/IP3/2016/74}

\subjclass[2010]{41A60, 33E17, 35Q15}

\keywords{Painlev\'e II equation, Riemann-Hilbert problem, modified Korteweg-de\thinspace Vries equation, asymptotic expansion}

\begin{abstract}
We are concerned with the defocusing modified Korteweg-de\thinspace Vries equation equipped with particular type of irregular initial conditions that are given as linear combinations of the Dirac delta function and Cauchy principal value. 
For the initial value problem we prove the existence of smooth self-similar solution, whose profile function is the Ablowitz-Segur solution of  the second Painlev\'e equation. Our method is to use the Riemann-Hilbert approach to improve asymptotics of these Painlev\'e II transcendents and find desired profile function by constructing its Stokes multipliers. 
\end{abstract}

\maketitle

\section{Introduction}
In this paper we study the Cauchy problem for the defocusing modified Korteweg-de\thinspace Vries (mKdV) equation
\begin{equation}\label{mkdv}
u_t+u_{xxx}-\frac{3}{2}u^{2}u_{x}=0,\quad  t>0, \ x\in\R,
\end{equation}
equipped with the initial conditions of the following particular form 
\begin{align}\label{init-cond}
u_{a,b}(x) := a\hspace{0.4pt}\delta(x)+b\hspace{1.5pt}\mathrm{p.v.}(1/x),\quad x\in\R,
\end{align}
where the parameters $a$, $b$ are real numbers, $\delta(x)$ denotes the Dirac delta function and $\mathrm{p.v.}\,(1/x)$ is the Cauchy principal value. 
Straightforward calculations show that, if $u$ is the solution of the equation \eqref{mkdv}, then the function 
\begin{align}\label{scale2}
u_{\lambda}(t,x) := \lambda u(\lambda^{3}t,\lambda x),\quad\lambda>0
\end{align} 
also has this property and consequently, the scaling argument (see e.g. \cite{ponce}, \cite{tao}) suggests the existence of solutions of the mKdV equation in the classical Sobolev space $H^{s}(\R)$ for $s\ge -1/2$. In the recent years, many effort has been made to develop the well-posedness theory for the equation \eqref{mkdv}. In particular, in the paper \cite{MR1211741}, the local existence, uniqueness and $C^{0}$-uniform continuity of solutions in the terms of the initial data was established in the space $H^{s}(\R)$ for $s\ge 1/4$.
The exponent $s=1/4$ appears to be optimal due to \cite{MR1813239}, where the ill-posedness (in the $C^{0}$-uniform sense) was showed for $s<1/4$. On the other hand, in \cite{MR1683054}, the results concerning global well-posedness of the equation \eqref{mkdv} were provided in the space $H^{s}(\R)$ for $s>3/5$. These in turn were improved in \cite{MR1969209} to the case $s>1/4$, by the application of the Miura transform and $I$-method for almost conservation laws. Finally, the same techniques were applied in \cite{MR2531556} to establish global well-posedness for the exponent $s=1/4$. 
An alternative scale of function spaces for studying the existence of solutions for the mKdV equation, were introduced in \cite{MR2096258} by the norm
\begin{align}\label{norm1}
\|u\|_{\widehat{H}^{r}_{\sigma}(\R)} := \|\<\xi\>^{\sigma}\widehat{u}(\xi)\|_{L_{\xi}^{r'}(\R)},
\end{align}
where $r\ge 1$, $\sigma\ge 0$ and $1/r+1/r'=1$.
Then the combination of the results from \cite{MR2096258} and \cite{MR2529909} provides the local well-posedness (in locally Lipschitz sense) for the equation \eqref{mkdv} in the space $\widehat{H}^{r}_{\sigma}(\R)$, where $r\in(1,2)$ and $\sigma=\sigma(r):=1/2 - 1/(2r)$. 
It is known that the borderline pair $(r,\sigma) = (1,0)$ corresponds to the space, which is critical with respect to the scaling transformation \eqref{scale2} and, to the best of our knowledge, the well-posedness theory remains an open question in this case.
In this paper we are interested in the initial conditions of the form \eqref{init-cond} that are particular type of critical initial data from the space $\widehat{H}^{1}_{0}(\R)$. It is well-known that, given $\alpha\in\R$, meromorphic functions satisfying the second Painlev\'e (PII) equation 
\begin{align}\label{PII}
v_{xx}(x) = x v(x) + 2v^{3}(x) - \alpha,\quad x\in\C,
\end{align}
provide us an important class of self-similar solutions for the equation \eqref{mkdv}. To be more precise, if $v(x)$ is the PII transcendent which is pole-free on the real line and satisfies $v(x)\in\R$ for all $x\in\R$, then the function 
\begin{align}\label{self-sim}
u(t,x) := -2(3t)^{-1/3}v(x(3t)^{-1/3}),\quad t>0, \ x\in\R,
\end{align}
is a real-valued solution of the mKdV equation (see e.g. \cite{MR0481656}, \cite{MR2501035}, \cite{MR1207209}). The main result of this paper is the following theorem. 
\begin{theorem}\label{th-kdv}
Given $a\in\R$ and $b\in(-1,1)$, there is a solution $v$ of the second Painlev\'e equation such that $v$ is pole-free on the real line, $v(x)\in\R$ for $x\in\R$ and the corresponding function \eqref{self-sim} is a smooth solution of the equation \eqref{mkdv} satisfying the initial condition
\begin{align}\label{conv-1}
\lim_{t\to 0^{+}}u(t,x)= a\hspace{0.4pt}\delta(x)+b\hspace{1.5pt}\mathrm{p.v.}\,(1/x)\quad \text{in}\quad \mathcal{S}'(\R).
\end{align}
Furthermore, the above limit becomes pointwise in the frequency space, that is,
\begin{equation}\label{eq-lim-2}
\lim_{t\to 0^{+}}\widehat{u}(t,\xi) = a - i\pi b\,\mathrm{sgn}\,\xi,\quad \xi\in\R\setminus\{0\}.
\end{equation}
\end{theorem}
In the proof of the above theorem we use the approach based on the Riemann-Hilbert (RH) problem associated with the PII equation, that was originated in the fundamental papers \cite{MR0588248} and \cite{jimbo}. To be more precise, let us assume that $\alpha\in(-1/2,1/2)$ is a given number and let $\Sigma$ be a contour in the complex $\lambda$-plane, consisting of the six rays that are oriented from the origin to the infinity
$$\gamma_{k}:=\{\lambda\in\C\setminus\{0\} \ | \ \mathrm{arg}\, \lambda = \pi/6 + (k-1)\pi/3\},\quad k=1,2,\ldots,6.$$ 
The RH problem on the graph $\Sigma$ is defined by the {\em Stokes multipliers}, that is, the triple of complex numbers $(s_{1},s_{2},s_{3})$ satisfying the following constraint condition 
\begin{equation}\label{stokes2}
s_{1}-s_{2}+s_{3} + s_{1}s_{2}s_{3} = -2\sin(\pi\alpha).
\end{equation}
Then, any choice of the Stokes initial data, gives us a solution $\Phi(\lambda,x)$ of the corresponding RH problem, which is a $2\times 2$ matrix valued mapping, sectionally holomorphic in $\lambda$ and meromorphic with respect to the variable  $x$ (see \cite{MR0588248}, \cite{MR2264522}, \cite{jimbo} for more details). If we write $\theta(\lambda,x) := i(4\lambda^{3}/3 + x\lambda)$ and assume that $\sigma_{3}:=\mathrm{diag}\,(1,-1)$ is the third Pauli matrix, then the function $v(x)$ given by the limit 
\begin{align*}
v(x)=\lim_{\lambda\to\infty}(2\lambda\Phi(\lambda,x)e^{\theta(\lambda,x)\sigma_{3}})_{12},
\end{align*}
is a solution of the PII equation \eqref{PII}. Thus, we can define the mapping
\begin{align*}
\{(s_{1},s_{2},s_{3})\in\C^{3} \text{ satisfying }\eqref{stokes2}\}\to \{\text{solutions of the equation \eqref{PII}}\},
\end{align*}
which appears to be bijection between the set of all Stokes multipliers and the set of the Painlev\'e II transcendents (see e.g. \cite{MR2264522}). 
In this paper, we are interested in the real Ablowitz-Segur solutions corresponding to the following Stokes initial data
\begin{gather}\label{stokes-11bb}
s_1 = -\sin(\pi\alpha) -i k, \quad s_2=0,\quad s_3 = -\sin(\pi\alpha) + i k,\\ \label{stokes-11bbc}
\alpha\in(-1/2,1/2),\quad k\in(-\cos(\pi\alpha),\cos(\pi\alpha)),
\end{gather}
that, for the brevity, are denoted by $v(x;\alpha,k)$. It is well-known that the solutions are such that $v(x;\alpha,k)\in\R$ \ for \ $x\in\R$ (see \cite[Chapter 11]{MR2264522}) and, by the result \cite[Theorem 2]{MR3670014}, they are pole-free on the real axis. 
Furthermore, we have the following asymptotic behaviors
\begin{gather}\label{asy-f-plus}
v(x;\alpha,k)=\alpha x^{-1}+2\alpha(1-\alpha^2)x^{-4}+O(x^{-7}),\quad x\to\infty,\\ \label{asy-f-real}
\hspace{-4pt}v(x;\alpha,k) \!=\! \frac{d}{(-x)^{\frac{1}{4}}}\cos(\frac{2}{3}(-x)^{\frac{3}{2}} \!-\! \frac{3}{4}d^{2}\ln(-x)\!+\!\phi) \!+\! O((-x)^{-1}), \ \ x\to -\infty, 
\end{gather}
where the constants $d$ and $\phi$ representing the magnitude and phase shift of the leading term in \eqref{asy-f-real}, respectively, are given by the following {\em connection formulas}
\begin{gather}\label{conn-f-real-1}
d(k,\alpha) := \frac{1}{\sqrt{\pi}}\sqrt{-\ln(\cos^{2}(\pi\alpha) - k^{2})},\\ \label{conn-f-real-2}
\phi(k,\alpha) := -\frac{3}{2}d^{2}\ln2 + \mathrm{arg}\,\Gamma\left(\frac{1}{2}id^{2}\right) - \frac{\pi}{4} - \mathrm{arg}\,(-\sin(\pi\alpha) - ki).
\end{gather} 
The above asymptotic relations and connection formulas were formally obtained in \cite{MR0859353} and rigorously justified in \cite{MR1174016} by the isomonodromy method. Another rigorous proof of \eqref{asy-f-plus}, based on the steepest descent analysis of the RH problem associated with the PII equation, were provided in \cite{MR1950792}. The same techniques were also successfully applied in \cite{MR3670014} to establish the asymptotic \eqref{asy-f-real} together with the connection formulas \eqref{conn-f-real-1} and \eqref{conn-f-real-2}. 
In the proof of Theorem \ref{th-kdv}, we provide explicit formulas for the parameters $\alpha$ and $k$, in the terms of the coefficients $a\in\R$ and $b\in(-1,1)$, such that the corresponding function $v(x;\alpha,k)$ is the expected profile function for the self-similar solution of the equation \eqref{mkdv}, satisfying the initial condition \eqref{conv-1}. To this end, we will need the following result, which develops the remainder term $O((-x)^{-1})$ from the asymptotic relation \eqref{asy-f-real}.
\begin{theorem}\label{th-asymptotic-real}
Given $\alpha\in (-1/2,1/2)$ and $k\in (-\cos(\pi\alpha), \cos(\pi\alpha))$, the corresponding real Ablowitz-Segur solution $v(x;\alpha,k)$ of the second Painlev\'e equation has the following asymptotic behavior as $x\to-\infty$:
\begin{align}\label{asym-1a}
v(x;\alpha,k) = \frac{d}{(-x)^{\frac{1}{4}}}\cos(\frac{2}{3}(-x)^{\frac{3}{2}} - \frac{3}{4}d^{2}\ln(-x)+\phi) +\frac{\alpha}{x} + O((-x)^{-\frac{7}{4}}),
\end{align}
where the constants $d$ and $\phi$ are given by the connection formulas \eqref{conn-f-real-1} and \eqref{conn-f-real-2}.
\end{theorem}
In the proof of the above theorem, we change the variables of the RH problem associated with the equation \eqref{PII} and recall transformations from \cite[p.16-19]{MR3670014}, leading to a RH problem, which is suitable for the use of the steepest descent techniques. In the new coordinates the phase function has the form $\tilde\theta(z):=i(4z^{3}/3 - z)$ and admits two critical points $z_{\pm}:=\pm 1/2$. Then, the contribution to the formula \eqref{asym-1a} coming from the part of the graph of the transformed RH problem, located away from the stationary points and the origin, is exponentially small. It is also known  that the local parametrices in neighborhoods of the points $z_{\pm}$ can be constructed explicitly using the parabolic cylinder functions. Consequently, in \cite{MR3670014} (see also \cite[p.\,322]{MR2264522}), it was shown that the leading term of \eqref{asym-1a} together with the connection formulas \eqref{conn-f-real-1}, \eqref{conn-f-real-2} can be derived from the asymptotic behavior of these special functions. This in turn gives precisely the relation \eqref{asy-f-real}.
The crucial point of the proof of Theorem \ref{th-asymptotic-real} is to improve \eqref{asy-f-real} and obtain the $\alpha x^{-1}$ term of the asymptotic \eqref{asym-1a}, using the explicit form of the parametrix near the origin, that was established in \cite[Theorem 6.5]{kok}. 
Then Theorem \ref{th-kdv} will be derived by the application of the relation \eqref{asym-1a} and {\em total integral formula} for the real Ablowitz-Segur solution $v(x;\alpha,k)$, which express the Cauchy principal value integral of the PII transcendent in the terms of the parameters $\alpha$ and $k$ (see Theorem \ref{th-total}).\\[5pt]
\noindent {\bf Outline.} The paper is organized as follows. In Section 2 we formulate the RH problem for the PII equation and recall its transformations leading to a RH problem, which is suitable for the use of the steepest descent techniques. In Section 3, we establish estimates for the solution of the transformed RH problem and provide a representation of its solution in the terms of appropriate local parametrices. Section 4 is devoted for the proof of Theorem \ref{th-asymptotic-real}, whereas in Section 5 we prove Theorem \ref{th-kdv}. The last section is the Appendix, where we consider the local parametrices that are required in the analysis of the transformed RH problem.\\[5pt]
{\bf Notation and terminology.} We define $M_{2\times 2}(\C)$ to be the complex linear space consisting of the $2\times 2$ complex matrices, endowed with the Frobenius norm 
\begin{align*}
\|A\|:=\sqrt{|a_{11}|^{2} + |a_{12}|^{2} + |a_{21}|^{2} + |a_{22}|^{2}},\quad A=[a_{lm}]\in M_{2\times 2}(\C).
\end{align*}
It is known that the norm is sub-multiplicative, that is, 
\begin{equation}\label{from-mult}
\|AB\|\le \|A\| \|B\|,\quad A,B\in M_{2\times 2}(\C).
\end{equation}
If $\Sigma$ is a contour contained in the complex plane and $1\le p <\infty$, then $L^{p}(\Sigma)$ is the space of measurable functions $f:\Sigma\to M_{2\times 2}(\C)$, equipped with the usual norm
\begin{align*}
\|f\|_{L^{p}(\Sigma)}:= \left(\int_{\Sigma}\|f(z)\|^{p}\,|dz|\right)^{1/p}.
\end{align*}
Furthermore, for $p=\infty$ the norm takes the following form 
\begin{align*}
\|f\|_{L^{\infty}(\Sigma)}:= \mathrm{ess\,sup}_{z\in \Sigma}\,\|f(z)\|.
\end{align*}
If $1\le p<\infty$ and the contour $\Sigma$ is unbounded, then we follow the notation from \cite{zhou} and define the space $L^{p}_{I}(\Sigma)$ consisting of functions $f:\Sigma\to M_{2\times 2}(\C)$ with the property that there is $f(\infty)\in M_{2\times 2}(\C)$ such that $f-f(\infty)\in L^{p}(\Sigma)$. It is clear that the matrix $f(\infty)$ is uniquely determined by $f$ and hence, we can set norm 
\begin{align*}
\|f\|_{L^{p}_{I}(\Sigma)}:=\left(\|f-f(\infty)\|^{p}_{L^{p}(\Sigma)} + \|f(\infty)\|^{p}\right)^{1/p},\quad f\in L^{p}_{I}(\Sigma).
\end{align*}
Throughout this paper we frequently write $a\lesssim b$ to denote the inequality $a\le cb$ for some $c>0$. Furthermore, we use the notation $a\sim b$ provided there are constants $c_{1},c_{2}>0$ such that $c_{1}b \le a\le c_{2}b$.

\section{The RH approach for the PII equation}

\begin{figure}[h]
\begin{center}
\includegraphics[scale=0.62]{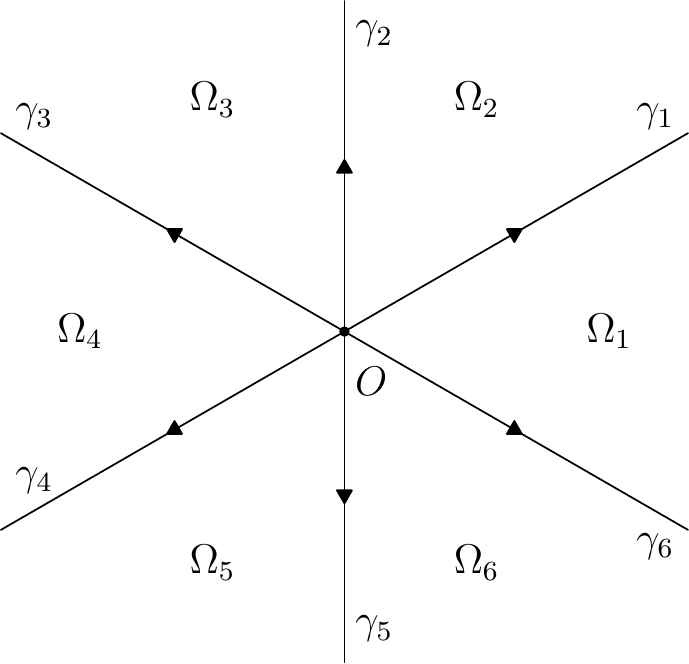}
\end{center}
\caption{Contour of the RH problem associated with the inhomogeneous PII equation.}
\label{fig1}
\end{figure}

\begin{figure}[h]
\centering
\includegraphics[scale=0.65]{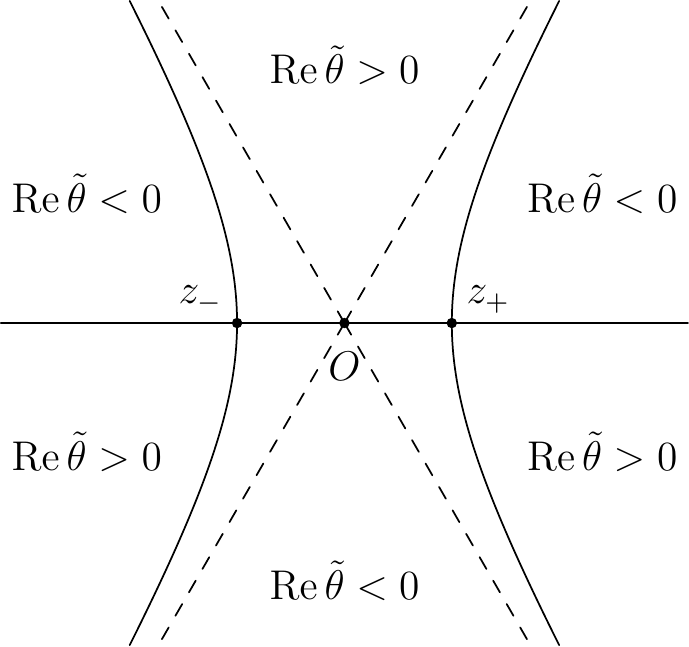}
\caption{The regions of sign changing of the function $\mathrm{Re}\,\tilde\theta(z)$. The dashed rays have directions $\exp(ik\pi/3)$ for $k=1,2,4,5$.}
\label{fig2}
\end{figure}

\begin{figure}[h]
\begin{subfigure}{0.49\textwidth}
\begin{center}
\includegraphics[scale=0.63]{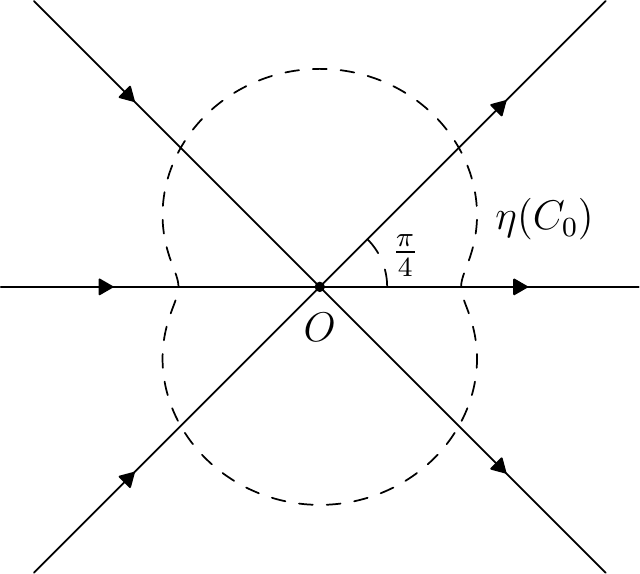}
\end{center}
\end{subfigure}
\begin{subfigure}{0.49\textwidth}
\begin{center}
\includegraphics[scale=0.63]{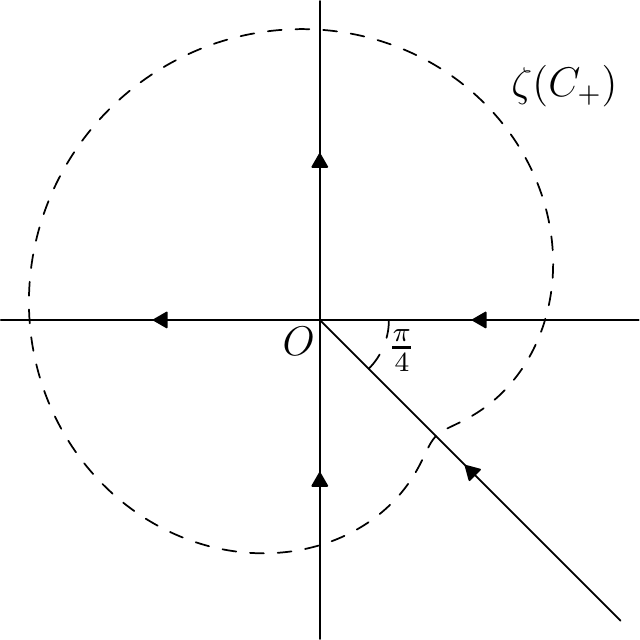}
\end{center}
\end{subfigure}
\caption{Left: the contour $\Sigma^{0}_{T}$ and the closed curve $\eta(C_{0})$. 
Right: the contour $\Sigma^{+}_{T}$ with the closed curve $\zeta(C_{+})$.}
\label{fig3}
\end{figure}

\begin{figure}[h]
\begin{center}
\includegraphics[scale=1.11]{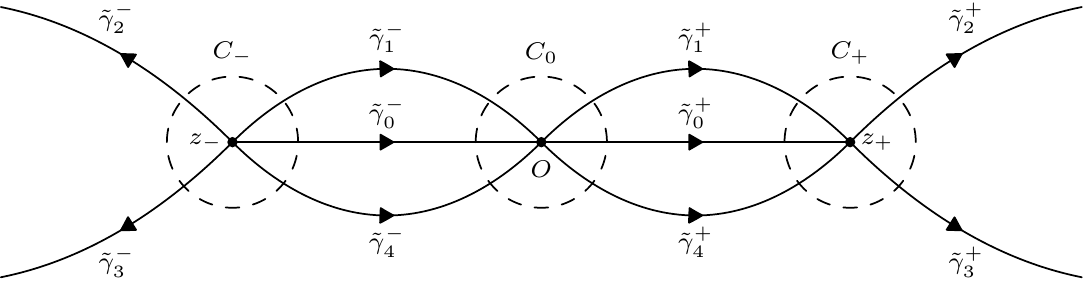}
\end{center}
\caption{The contour $\Sigma_{T}$ and the circles $C_{0}$, $C_{\pm}$ that are depicted by dashed lines.}
\label{fig4}
\end{figure}

\begin{figure}[h]
\centering
\includegraphics[scale=1.1]{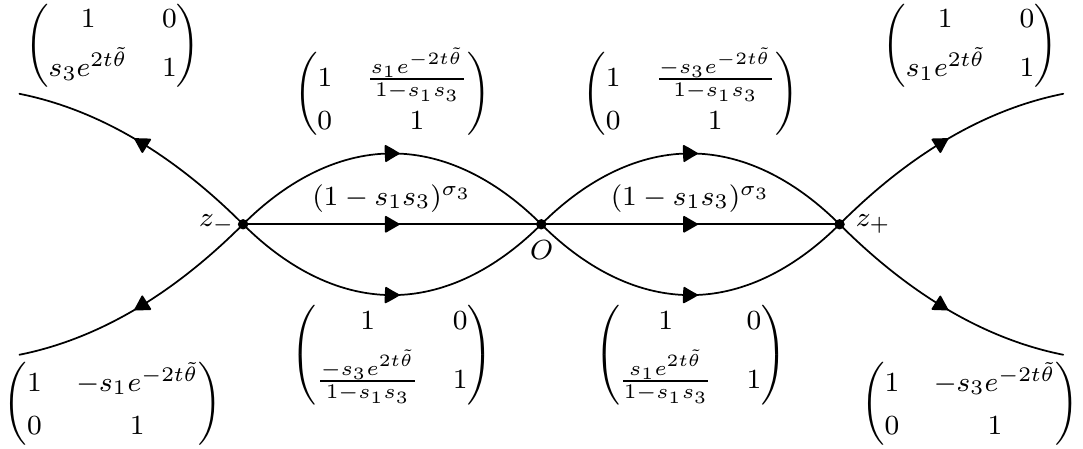}
\caption{The contour $\Sigma_{T}$ and the associated jump matrices for the RH problem fulfilled by the function $T(z)$.}
\label{fig5}
\end{figure}

In this section we intend to formulate the Riemann-Hilbert problem for the inhomogeneous PII equation \eqref{PII} and recall an approach, relying on 
normalization of jump matrices at infinity and contour deformation, that will allow us to apply the steepest descent techniques. To this end, let us consider the contour $\Sigma$ in the complex $\lambda$-plane consisting of the six rays $$\gamma_{j}:\quad \mathrm{arg}\, \lambda = \pi/6 + (j-1)\pi/3, \qquad j=1,2,\ldots,6,$$ that are oriented from zero to infinity, as it is depicted on Figure \ref{fig1}. The contour divides the complex plane on the six regions that we denote by $\Omega_{1},\Omega_{2},\ldots,\Omega_{6}$. Furthermore, due to the orientation, we can easily distinguish the left $(+)$ and right $(-)$ sides of the graph $\Sigma$. For any $1\le j\le 6$, each of the rays $\gamma_{j}$ has assigned a triangular {\em jump matrix} $S_{j}$, given by $$S_{j}:=\begin{pmatrix}1& 0\\ s_{j} & 1\end{pmatrix}, \quad j=1,3,5 \quad\text{and}\quad S_{j}:=\begin{pmatrix}1& s_{j}\\ 0 & 1\end{pmatrix}, \quad j=2,4,6,$$ where the parameters $(s_{1}, s_{2},\ldots,s_{6})\in\C^{6}$ are called the {\em Stokes multipliers} and satisfy $s_{k} = s_{k+3}$ for $k=1,2,3$ together with the constraint condition \eqref{stokes2}. Let us assume that $\sigma_1$, $\sigma_2$ and $\sigma_3$ denote the Pauli matrices given by
\begin{align*}
\sigma_{1}:=\begin{pmatrix} 0& 1\\ 1& 0 \end{pmatrix},\quad \sigma_{2}:=\begin{pmatrix}0& -i\\ i& 0 \end{pmatrix},\quad
\sigma_{3}:=\begin{pmatrix}1& 0\\ 0&-1\end{pmatrix}. 
\end{align*}
Given $\alpha\in(-\tfrac{1}{2},\tfrac{1}{2})$, the Riemann-Hilbert problem associated with the PII equation consists in finding a function $\Phi(\lambda)=\Phi(\lambda;x)$ with the values in the matrix space $M_{2\times 2}(\C)$ such that the following conditions are satisfied.\\[5pt]
\noindent\makebox[5.5mm][l]{$(1)$}\parbox[t][][t]{122mm}{The function $\Phi(\lambda)$ is analytic for $\lambda\in\C\setminus\Sigma$.}\\[4pt]
\noindent\makebox[5.5mm][l]{$(2)$}\parbox[t][][t]{122mm}{Given $\lambda\in\Sigma\setminus\{0\}$, let $\Phi_{+}(\lambda)$ and $\Phi_{-}(\lambda)$ denote the limiting values of $\Phi(\lambda')$ as $\lambda'$ tends to $\lambda$ from the left and right side of the contour $\Sigma$, respectively. Then, for any $1\le j\le 6$, we have the following jump relation $$\Phi_{+}(\lambda) = \Phi_{-}(\lambda) S_{j}, \quad\lambda\in\gamma_{j}.$$}\\[4pt]
\noindent\makebox[5.5mm][l]{$(3)$}\parbox[t][][t]{122mm}{The function $\Phi(\lambda)$ has the following asymptotic behavior
$$\Phi(\lambda) = (I + O(\lambda^{-1}))e^{-\theta(\lambda)\sigma_{3}}, \quad \lambda\to\infty,$$ where $\theta(\lambda,x) := i(4\lambda^{3}/3 + x\lambda)$ is a phase function.}\\[4pt]
\noindent\makebox[5.5mm][l]{$(4)$}\parbox[t][][t]{122mm}{We have the following asymptotic relation 
\begin{align*}
\Phi(\lambda)= O\begin{pmatrix}|\lambda|^{-|\alpha|}& |\lambda|^{-|\alpha|}\\[5pt] |\lambda|^{-|\alpha|}& |\lambda|^{-|\alpha|}\end{pmatrix},\quad \lambda\to 0.
\end{align*}}\\[5pt]
Throughout this paper we consider the real Ablowitz-Segur solutions corresponding to the Stokes initial data given by \eqref{stokes-11bb} and \eqref{stokes-11bbc}. This in turn implies that $S_{2} = S_{5} = I$. Changing the variables on the complex $\lambda$-plane by the formulas
\begin{equation}\label{change-var}
\lambda(z) = (-x)^{1/2}z, \ \ t(x)=(-x)^{3/2},\quad z\in \C, \ x<0,
\end{equation}
we obtain $\theta(\lambda,x) = t\tilde \theta(z)$, where $\tilde\theta(z) := i(4z^{3}/3 - z)$. Let us assume that $\tilde\Sigma$ is a contour in the complex $z$-plane, consisting of the four oriented rays $$\gamma_{k}:\quad\mathrm{arg}\, \lambda = \pi/6 + (k-1)\pi/3, \qquad k=1,3,4,6,$$
that are equipped with the following triangular jump matrices 
\begin{align*}
G_{2k}:= \begin{pmatrix} 1 & e^{-2t\tilde\theta(z)}s_{2k}\\[5pt] 0 & 1\end{pmatrix}, \  k=2,3 \ \ \text{and} \ \ 
G_{2k+1}:=\begin{pmatrix} 1 & 0\\[5pt] e^{2t\tilde\theta(z)}s_{2k+1} & 1\end{pmatrix}, \  k=0,1.
\end{align*}
Then, it is not difficult to check that the function 
\begin{align*}
U(t,z) :=\Phi(\lambda(z),-t^{2/3})\exp(t\tilde\theta(z)\sigma_{3})
\end{align*}
is a solution of the following normalized Riemann-Hilbert problem. \\[3pt]
\noindent\makebox[8.5mm][l]{$(D1)$}\parbox[t][][t]{119mm}{The function $U(z)$ is holomorphic for $z\in \C\setminus\tilde\Sigma$.}\\[2pt]
\noindent\makebox[8.5mm][l]{$(D2)$}\parbox[t][][t]{119mm}{For any $k=1,3,4,6$, we have the following jump condition $$U_{+}(z) = U_{-}(z)G_{k},\quad z\in \gamma_{k}.$$}\\
\noindent\makebox[8.5mm][l]{$(D3)$}\parbox[t][][t]{119mm}{The function $U(z)$ satisfies the following asymptotic relation 
$$U(z) = I + O(z^{-1}),\quad z\to\infty.$$}\\
\noindent\makebox[8.5mm][l]{$(D4)$}\parbox[t][][t]{122mm}{We have the following asymptotic relation 
\begin{align*}
U(z)= O\begin{pmatrix}|z|^{-|\alpha|}& |z|^{-|\alpha|}\\[5pt] |z|^{-|\alpha|}& |z|^{-|\alpha|}\end{pmatrix},\quad z\to 0.
\end{align*}}\\[5pt]
Let us observe that the scaled phase function $\tilde\theta(z)$ has two stationary points $z_{\pm} := \pm 1/2$ such that $\tilde\theta(\pm 1/2)=\mp i/3$. 
Therefore the real line and the curves $$h_{\pm}(t):=it \pm\sqrt{t^{2}/3 +1/4},\quad t\in\R$$ are solutions of the equation $\mathrm{Re}\,\tilde\theta(z)=0$ passing through the stationary points $z_{\pm}$. Clearly the curves $h_{+}$ and $h_{-}$ are asymptotic to the rays $\mathrm{arg}\,\lambda = \pm\frac{\pi}{3}$ and $\mathrm{arg}\,\lambda = \pm\frac{2\pi}{3}$, respectively, and together with the real axis they separate the regions of the sign changing of the function $\mathrm{Re}\,\tilde\theta(z)$, as it is depicted on Figure \ref{fig2}.
By the results of \cite[Chapter 3]{MR3670014}, we can use the sign changing regions of the function $\mathrm{Re}\,\tilde\theta(z)$ and define the equivalent RH problem on a contour $\Sigma_{T}$.
To describe the contour more precisely we will use two auxiliary graphs $\Sigma_{T}^{0}$ and $\Sigma_{T}^{+}$ that are depicted on Figure \ref{fig3}. 
The former graph consists of the six rays 
\begin{align*}
\mathrm{arg}\,\lambda = 0, \quad \mathrm{arg}\,\lambda = \pi, \quad \mathrm{arg}\,\lambda = \pi/4 + j\pi/2,\quad 0\le j\le 3
\end{align*}
and the later is formed by the curves
\begin{align*}
\mathrm{arg}\,\lambda = 7\pi/4, \quad \mathrm{arg}\,\lambda = j\pi/2,\quad 0\le j\le 3.
\end{align*}
For the local change of coordinates, we will also need the mappings $\eta(z)$ and $\zeta(z)$ that are defined, in a neighborhood of the origin and the point $z_{+}$, respectively, by the following formulas 
\begin{gather}\label{d-zeta}
\eta(z):=i\tilde\theta(z) = z - 4z^3/3,\\ \label{d-eta}
\zeta(z):= 2\sqrt{-\tilde\theta(z)+\tilde\theta(z_{+})}=4\sqrt{3}e^{\frac{3}{4}\pi i}\left(z-1/2\right)(z+1)^{\frac{1}{2}}/3,
\end{gather}
where the branch cut of the square root is taken such that $\mathrm{arg}\,(z-1/2)\in(-\pi,\pi)$. The functions $\eta(z)$ and $\zeta(z)$ are holomorphic in a neighborhood of the origin and $z_{+}$, respectively. Since $\eta'(0)\neq 0$ and $\zeta'(z_{+})\neq 0$, by the inverse mapping theorem, there is a sufficiently small $\delta>0$ \label{page-11} such that the functions $\eta(z)$ and $\zeta(z)$ are biholomorphic on the open balls $B(0,2\delta)$ and $B(z_{+},2\delta)$, respectively. \label{holo-eta} If we take $C_{0}:=\partial B(0,\delta)$ and $C_{\pm}:= \partial B(z_{\pm},\delta)$, then both $\eta(C_{0})$ and $\zeta(C_{+})$ are closed curves surrounding the origin (see Figure \ref{fig3}). We define $\Sigma_{T}$ to be a contour depicted on Figure \ref{fig4}, where $\t\gamma^{\pm}_{j}$, for $j=0,1,4$ are curves connecting the origin with the stationary points $z_{\pm}$ such that $\t\gamma^{\pm}_{0}$ are segments lying on the real line, while $\t\gamma^{\pm}_{0}$ and $\t\gamma^{\pm}_{4}$ are such that the sets $\t\gamma^{\pm}_{1}\setminus\{0,z^{\pm}\}$ and $\t\gamma^{\pm}_{4}\setminus\{0,z^{\pm}\}$ are contained in the lower and upper half-plane of $\C$, respectively. 
We also assume that $\t\gamma_{2}^{\pm}$ and $\t\gamma_3^{\pm}$ are unbounded components of the contour $\Sigma_{T}$ emanating from the stationary point $z_{\pm}$, that are asymptotic to the rays $\{\mathrm{arg}\, \lambda = \pi/2\mp\pi/3\}$ and $\{\mathrm{arg}\, \lambda = 3\pi/2 \pm \pi/3\}$, respectively. We require also that the part of $\Sigma_{T}$ contained in the ball $B(0,\delta)$ is the inverse image of the set $\Sigma^{0}\cap\eta(B(0,\delta))$ under the map $\eta$ restricted to the ball $B(0,2\delta)$ as well as the part of the contour $\Sigma_{T}$ contained in the ball $B(z_{+},\delta)$ is the inverse image of the set $\Sigma^{+}\cap\zeta(B(z_{+},\delta))$ under the map $\zeta$, restricted to the ball $B(z_{+},2\delta)$. Furthermore the part of the contour $\Sigma_{T}$ contained in the ball $B(z_{-},\delta)$ is taken to be a reflection across the origin of the set $\Sigma_{T}\cap B(z_{+},\delta)$. 
Then, by \cite[Section 3.1 and 3.2]{MR3670014}, the solution $U(t,z)$ can be deformed to the function $T(z):=T(t,z)$ with values in the space $M_{2\times 2}(\C)$, which satisfies the following RH problem on the graph $\Sigma_{T}$. \\[3pt]
\noindent\makebox[7.5mm][l]{$(T1)$} \parbox[t][][t]{118mm}{The function $T(z)$ is holomorphic for $z\in \C\setminus\Sigma_{T}$.}\\[3pt]
\noindent\makebox[7.5mm][l]{$(T2)$} \parbox[t][][t]{118mm}{Given $z\in\Sigma_{T}\setminus\{z_{+},z_{-},0\}$, let $T_{+}(z)$ and $T_{-}(z)$ denote the limiting values of $T(z')$ as $z'$ tends to $z$ from the left and right side of the contour $\Sigma_{T}$, respectively. Then, we have the jump relation 
$$T_{+}(z)= T_{-}(z)S_{T}(z),\quad z\in \Sigma_{T},$$
where the jump matrix $S_{T}(z)$ is presented on Figure \ref{fig5}.}\\[3pt]
\noindent\makebox[7.5mm][l]{$(T3)$} \parbox[t][][t]{118mm}{The function $T(z)$ has the following asymptotic behavior
$$T(z) = I + O(z^{-1}),\quad z\to\infty.$$}\\
\noindent\makebox[7.5mm][l]{$(T4)$} \parbox[t][][t]{118mm}{As $z\to z_{\pm}$, the function $T(z)$ is bounded.}\\[3pt]
\noindent\makebox[7.5mm][l]{$(T5)$} \parbox[t][][t]{122mm}{We have the following asymptotic relation 
\begin{align*}
T(z)= O\begin{pmatrix}|z|^{-|\alpha|}& |z|^{-|\alpha|}\\[5pt] |z|^{-|\alpha|}& |z|^{-|\alpha|}\end{pmatrix},\quad z\to 0.
\end{align*}}\\[5pt]
Observe that, from the construction of the function $T(z)$, it follows that the solution $v(x)$ of the PII equation \eqref{PII} can be obtained by the following limit
\begin{align}\label{u-lim}
v(x) = 2\sqrt{-x}\lim_{z\to\infty}(zT_{12}(z,(-x)^{-3/2})). 
\end{align}
\section{Representation of solutions of the RH problem on $\Sigma_{T}$} 
\begin{figure}[h]
\includegraphics[scale=0.8]{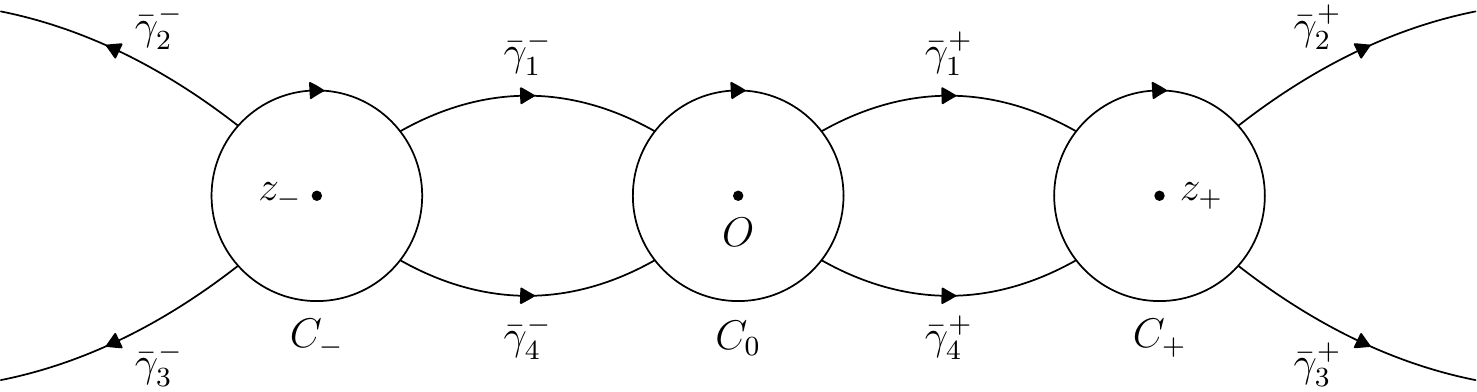}
\caption{The contour $\Sigma_{R}$ for the Riemann-Hilbert problem satisfied by the function $R(z)$.}
\label{fig6}
\end{figure}

In this section we provide representation for the solution of the RH problem (T1)--(T5), in the terms of local parametrices considered in the Appendix, and we establish estimates that will be used in the proof of Theorem \ref{th-asymptotic-real}. Let us assume that $\Sigma_{R}$ is the contour in the complex plane consisting of the circles $C_{\pm}=\partial D(z_{\pm},\delta)$, $C_{0}=\partial D(z_{0},\delta)$ (see page \pageref{page-11}) and the parts $\bar \gamma^{\pm}_{j}$ of the curves $\tilde\gamma^{\pm}_{j}$ lying outside the set $B(z_{+},\delta)\cup B(z_{-},\delta)\cup B(z_{0},\delta)$ (see Figure \ref{fig6}). Let $T^{(0)}(z)$ be the local parametrix near the origin, defined by the formula \eqref{equa-p}, and let $T^{(r)}(z)$ (resp. $T^{(l)}(z)$) be the parametrix near the stationary point $z_{+}=1/2$ (resp. $z_{-}=-1/2$), given by the equation \eqref{t-r-def} (resp. \eqref{t-l-def}). We consider the map $R(z)$ represented as follows
\begin{gather*}
R(z):=T(z)T^{(0)}(z)^{-1}, \ z\in D(0,\delta)\setminus\Sigma_{T},\\
R(z):=T(z)T^{(r)}(z)^{-1}, \, z\in D(z_{+},\delta)\setminus\Sigma_{T}, \ R(z):=T(z)T^{(l)}(z)^{-1}, \, z\in D(z_{-},\delta)\setminus\Sigma_{T},\\
R(z):=T(z)N(z)^{-1}, \ z\in\C\setminus(D(z_{+},\delta)\cup D(z_{-},\delta)\cup D(0,\delta)\cup\Sigma_{T}),
\end{gather*}
where $N(z)$ is the function defined by the formula \eqref{f-z}. In view of the equality \eqref{u-lim} and the fact that $N(z) = I+O(1/z)$ as $z\to\infty$, it is not difficult to check that the solution of the corresponding PII equation can be obtained by the limit
\begin{align}\label{u-lim-2}
v(x) = 2\sqrt{-x}\lim_{z\to\infty}(zR_{12}(z,(-x)^{-3/2})).
\end{align}
Let $S_{R}(z)$ be the jump matrix on the contour $\Sigma_{R}$, given by
\begin{gather*}
S_{R}(z):=T^{(0)}(z)N(z)^{-1}, \quad z\in \partial D(0,\delta), \\
S_{R}(z):=T^{(r)}(z)N(z)^{-1},  \ \  z\in \partial D(z_{+},\delta),  \quad  S_{R}(z):=T^{(l)}(z)N(z)^{-1},  \ \ z\in \partial D(z_{-},\delta),\\
S_{R}(z):=N(z)S_{T}(z)N(z)^{-1}, \ z\in\Sigma_{R}\setminus(\partial D(z_{+},\delta)\cup \partial D(z_{-},\delta)\cup \partial D(0,\delta)).
\end{gather*}
Then the function $R(z)$ is a solution of the following RH problem. \\[3pt]
\noindent\makebox[7.5mm][l]{$(R1)$} \parbox[t][][t]{121mm}{The function $R(z)$ is analytic in $\C\setminus\Sigma_R$.}\\[3pt]
\noindent\makebox[7.5mm][l]{$(R2)$} \parbox[t][][t]{121mm}{We have the following jump relation
\begin{align*}
R_{+}(z) = R_{-}(z)S_{R}(z),\quad z\in \Sigma_{R}.
\end{align*}}\\
\noindent\makebox[7.5mm][l]{$(R3)$} \parbox[t][][t]{121mm}{The function $R(z)$ has the following asymptotic behavior
\begin{align*}
R(z)=I + O(z^{-1}),\quad z\to\infty.
\end{align*}}\\
It is known (see e.g. \cite{MR2264522}, \cite{zhou}) that the solution of the RH problem (R1)--(R3) can be obtained using the following fixed point problem 
\begin{align}\label{fix-point}
\rho= I+\RR(\rho),
\end{align}
where $I$ is the identity matrix and $\RR:L^{2}_{I}(\Sigma_{R})\to L^{2}_{I}(\Sigma_{R})$ is a linear map given by
\begin{equation*}
\RR(\rho):= \CC_{-}(\rho(S_{R} - I)),\quad \rho\in L^{2}_{I}(\Sigma_{R}).
\end{equation*}
In the above formula $\mathcal{C}_-:L^{2}(\Sigma_{R})\to L^{2}(\Sigma_{R})$ represents the Cauchy operator
\begin{align*}
[\CC_- f](z) := \lim_{z'\to z}\frac{1}{2\pi i}\int_{\Sigma_{R}} \frac{f(\xi)}{\xi-z'}\,d\xi, \quad z\in\Sigma_{R},
\end{align*}
where the parameter $z'$ tends non-tangentially to $z$ from the $(-)$ side of the graph $\Sigma_{R}$. It is known (see e.g. \cite[Section 2.5.4]{MR3450072}) that $\CC_{-}$ is a bounded operator on the space $L^2(\Sigma_{R})$. In particular, for any measurable sets $P_{1},P_{2} \subset\Sigma_{R}$, we have
\begin{equation}
\begin{aligned}\label{ineq-c}
\|\CC_{-}(f\chi_{P_{1}})\|_{L^{2}(P_{2})} & \le \|\CC_{-}(f\chi_{P_{1}})\|_{L^{2}(\Sigma_{R})}\le \|\CC_{-}\|_{L^{2}(\Sigma_{R})} \|f\chi_{P_{1}}\|_{L^{2}(\Sigma_{R})} \\
& = \|\CC_{-}\|_{L^{2}(\Sigma_{R})} \|f\|_{L^{2}(P_{1})},\quad f\in L^{2}(\Sigma_{R}).
\end{aligned}
\end{equation}
If the function $\rho\in L^2_{I}(\Sigma_{R})$ satisfies the equation \eqref{fix-point}, then the integral 
\begin{align}\label{ll-1}
R(z)= I +\frac{1}{2\pi i}\int_{\Sigma_{R}} \frac{\rho(\xi)(S_{R}(\xi)-I)}{\xi-z}\,d\xi, \quad z\not\in\Sigma_{R},
\end{align}
represents the solution of the RH problem defined on the contour $\Sigma_{R}$ and satisfies 
\begin{equation}\label{ll-2}
R_{-}(z)=\rho(z),\quad z\in\Sigma_{R}.
\end{equation}
In the following lemmata we provide estimates on the functions $S_{R}(z)$ and $R_{-}(z)$. 
\begin{lemma}\label{lem-est-3}
Let us define $\Sigma'_{R}:=\Sigma_R\setminus (C_{+}\cup C_{-} \cup C_{0})$. Then we have the following asymptotic relations
\begin{align}\label{bbb12kk}
&\|S_{R}-I\|_{L^{2}(C_{0})} = O(t^{-1}),&&\hspace{-50pt} t\to\infty,\\ \label{bbb1}
&\|S_{R}-I\|_{(L^{2}\cap L^{\infty})(\Sigma_{R})} = O(t^{-1/2}),&&\hspace{-50pt} t\to\infty,\\ \label{bbb12}
&\|S_{R}-I\|_{L^{2}(\Sigma'_{R})} = O(e^{-ct}),&&\hspace{-50pt} t\to\infty,
\end{align}
where $c>0$ is a constant.
\end{lemma}
\proof Applying Propositions \ref{l-par} and \ref{l-par2}, we infer that 
\begin{align*}
\|S_{R} - I\|_{L^{\infty}(C_{+})} = \|T^{(r)}N^{-1} - I\|_{L^{\infty}(C_{+})} = O(t^{-1/2}), \quad t\to\infty,\\
\|S_{R} - I\|_{L^{\infty}(C_{-})} = \|T^{(l)}N^{-1} - I\|_{L^{\infty}(C_{-})} = O(t^{-1/2}), \quad t\to\infty.
\end{align*}
Furthermore, by Theorem \ref{lok-aprox}, we have the following asymptotic
\begin{equation}\label{eq-mm-nn-1}
\|S_{R} - I\|_{L^{\infty}(C_{0})} = \|T^{(0)}N^{-1} - I\|_{L^{\infty}(C_{0})} = O(t^{-1}), \quad t\to\infty.
\end{equation}
Therefore, there is $t_{0}>0$ such that 
\begin{equation}
\begin{aligned}\label{ee-m}
&\|S_{R} - I\|_{L^{\infty}(C_{\pm})} \lesssim t^{-1/2}, \quad \|S_{R} - I\|_{L^{\infty}(C_{0})} \lesssim t^{-1},\quad t\ge t_{0},
\end{aligned}
\end{equation}
which implies that, for any $t>t_{0}$, we have
\begin{align}\label{ee-4-gg}
\|S_{R} - I\|^{2}_{L^{2}(C_{0})} = \int_{C_{0}}\|S_{R}(z) - I\|^{2}\,|dz|\lesssim \int_{C_{0}}t^{-2}\,|dz|\sim t^{-2} 
\end{align}
and hence the asymptotic relation \eqref{bbb12kk} follows. Using \eqref{ee-m} once again we infer that, for any $t>t_{0}$, the following inequality holds
\begin{align}\label{ee-4}
\|S_{R} - I\|^{2}_{L^{2}(C_{\pm})} = \int_{C_{\pm}}\|S_{R}(z) - I\|^{2}\,|dz|\lesssim \int_{C_{\pm}}t^{-1}\,|dz|\sim t^{-1}.
\end{align}
By the definition of $N(z)$ there is a constant $C>0$ such that 
\begin{align*}
\|N(z)\|\le C\quad\text{and}\quad \|N(z)^{-1}\|\le C,\quad z\in \Sigma_{R}'.
\end{align*}
Hence, using the inequality \eqref{from-mult}, we obtain
\begin{equation}
\begin{aligned}\label{eq-r-t}
&\|S_{R}(z) - I\| = \|N(z)[S_{T}(z) - I]N(z)^{-1}\|\\
& \quad\le \|N(z)\|\|S_{T}(z) - I\|\|N(z)^{-1}\|\le C^{2}\|S_{T}(z) - I\|,\quad z\in \Sigma_{R}'.
\end{aligned}
\end{equation}
Since the curve $\bar\gamma^{+}_{2}$ is asymptotic to the ray $\{se^{i\pi/6} \ | \ s>0\}$, there is $a >0$ and a smooth function $h:[a,\infty)\to\R$
satisfying the following asymptotic condition
\begin{equation*}
h(s)/s\to \sqrt{3}/3,\quad s\to\infty,
\end{equation*}
such that the map $\bar\gamma^{+}_{2}(s):=s + ih(s)$ for $s\ge a$ is the parametrization of the curve $\bar\gamma^{+}_{2}$. Let us take a small $\ve_{0}>0$ such that 
\begin{equation}\label{eq-ineq-55}
4(\sqrt{3}/3+\ve_{0})^{3}/3- 4(\sqrt{3}/3-\ve_{0})<0
\end{equation}
and observe that, there is $a_{0}> a$ such that, for any $s>a_{0}$, we have
\begin{align*}
&\mathrm{Re}\,\tilde\theta(s+ih(s))= 4 h(s)^{3}/3- 4s^{2}h(s)+h(s)\\
&\quad\le (4(\sqrt{3}/3+\ve_{0})^{3}/3- 4(\sqrt{3}/3-\ve_{0}))s^{3} + (\sqrt{3}/3+\ve_{0})s. 
\end{align*}
Therefore, in view of \eqref{eq-ineq-55}, there is $a_{1}>a_{0}>0$ with the property that 
\begin{align}\label{eq-hh-22}
\mathrm{Re}\,\tilde\theta(\bar\gamma^{+}_{2}(s)) = \mathrm{Re}\,\tilde\theta(s+ ih(s))\le -s,\quad s\ge a_{1}.
\end{align}
By the diagram depicted on Figure \ref{fig2}, we obtain the existence of $c_{0}>0$ such that $\mathrm{Re}\,\tilde\theta(\bar\gamma^{+}_{2}(s))\le -c_{0} $ for $s\in[a,a_{0}]$. Therefore, if we take $m:=\min(c_{0},a_{1})$, then 
\begin{align}\label{eq-hh-22bb}
\mathrm{Re}\,\tilde\theta(\bar\gamma^{+}_{2}(s))=\mathrm{Re}\,\tilde\theta(s+ ih(s))\le -m,\quad s\ge a.
\end{align}
Taking into account the form of the jump matrix $S_{T}$ on the curve $\tilde\gamma_{2}^{+}$ (see Figure \ref{fig5}), we infer that 
\begin{align}\label{eq-norm-1}
\|S_{T}(\bar\gamma^{+}_{2}(s)) - I\| = |s_{1}|e^{2t\mathrm{Re}\,\tilde\theta(\bar\gamma^{+}_{2}(s))},\quad s>a, \ t>0.
\end{align}
Combining this equality with \eqref{eq-r-t} and \eqref{eq-hh-22bb}, yields
\begin{align}\label{kk-ll-mm}
\|I - S_{R}\|_{L^{\infty}(\bar\gamma_{2}^{+})}\lesssim \|I - S_{T}\|_{L^{\infty}(\bar\gamma_{2}^{+})}\le \sup_{s\ge a} e^{2t\mathrm{Re}\,\tilde\theta(\bar \gamma^{+}_{2}(s))}\le e^{-2mt},\quad t>0.
\end{align}
Furthermore, using \eqref{eq-r-t}, \eqref{eq-hh-22}, \eqref{eq-hh-22bb} and \eqref{eq-norm-1}, we have 
\begin{equation}
\begin{aligned}\label{kk-ll-mm-22}
&\|I - S_{R}\|^{2}_{L^{2}(\bar\gamma_{2}^{+})}\lesssim \|I - S_{T}\|^{2}_{L^{2}(\bar\gamma_{2}^{+})}\sim\int_{a}^{\infty}|e^{2t\tilde\theta(\bar\gamma^{+}_{2}(s))}|^{2}|(\bar\gamma^{+}_{2})'(s)|\,ds\\
& \ \ \lesssim \int_{a}^{a_{1}}e^{4t\mathrm{Re}\,\tilde\theta(\bar\gamma_{2}^{+}(s))}\,ds + \int_{a_{1}}^{\infty}e^{4t\mathrm{Re}\,\tilde\theta(\bar\gamma_{2}^{+}(s))}\,ds \\
& \ \ \le \int_{a}^{a_{1}}e^{-4mt}\,ds+\int_{a_{1}}^{\infty}e^{-4ts}\,ds = (a_{1}-a)e^{-4mt} + (4t)^{-1}e^{-4ta_{1}}, \ t>0.
\end{aligned}
\end{equation}
Then \eqref{kk-ll-mm} and \eqref{kk-ll-mm-22} imply that there is a constant $c_{1}>0$ such that 
\begin{align}\label{asymp-s-r-1}
\|I - S_{R}\|_{(L^{2}\cap L^{\infty})(\bar\gamma_{2}^{+})} = O(e^{-c_{1}t}),\quad t\to\infty.
\end{align}
If $\bar\gamma^{+}_{1}:[0,1]\to\C$ is the parametrization of the curve $\bar\gamma^{+}_{1}$, then using the diagram from Figure \ref{fig2} once again, we obtain the existence of $m_{1}>0$ such that 
\begin{align}\label{est-kk-tt}
\mathrm{Re}\,\tilde\theta(\bar\gamma^{+}_{1}(s))\ge m_{1},\quad s\in[0,1].
\end{align}
Using the form of the jump matrix $S_{T}$ on the curve $\tilde\gamma_{1}^{+}$ (see Figure \ref{fig5}), we obtain
\begin{align}\label{eq-norm-2}
\|S_{T}(\bar\gamma^{+}_{1}(s)) - I\| = \frac{|s_{3}|}{|1-s_{1}s_{3}|}e^{-2t\mathrm{Re}\,\tilde\theta(\bar\gamma^{+}_{1}(s))},\quad 
s\in[0,1], \ t>0.
\end{align}
Combining this equality with \eqref{eq-r-t} and \eqref{est-kk-tt} gives
\begin{align}\label{ee-3}
\|I - S_{R}\|_{L^{\infty}(\bar\gamma_{1}^{+})} \lesssim \|I - S_{T}\|_{L^{\infty}(\bar\gamma_{1}^{+})} \le \sup_{s\in[0,1]} e^{-2t\mathrm{Re}\,\tilde\theta(\bar\gamma^{+}_{1}(s))}\le e^{-2m_1 t}, \  t>0
\end{align}
and furthermore, by \eqref{eq-r-t}, \eqref{est-kk-tt} and \eqref{eq-norm-2}, we have
\begin{equation}
\begin{aligned}\label{ee-2}
&\|I - S_{R}\|^{2}_{L^{2}(\bar\gamma_{1}^{+})}\lesssim \|I - S_{T}\|^{2}_{L^{2}(\bar\gamma_{1}^{+})}\sim \int_{0}^{1}e^{-4t\tilde\theta(\bar\gamma^{+}_{1}(s))}|\bar\gamma^{+}_{1}(s)'|\,ds\\
&\qquad \lesssim \int_{0}^{1}e^{-4t\tilde\theta(\bar\gamma^{+}_{1}(s))}\,ds \le e^{-4m_1t},\qquad t>0.
\end{aligned}
\end{equation}
Combining \eqref{ee-3} and \eqref{ee-2}, we deduce that
\begin{align}\label{asymp-s-r-2}
\|I - S_{R}\|_{(L^{2}\cap L^{\infty})(\bar\gamma_{1}^{+})} = O(e^{-c_{2}t}),\quad t\to\infty,
\end{align}
where $c_{2}:=2m_{1}$. Let us observe that we can repeat the above argument to obtain the asymptotics \eqref{asymp-s-r-1} and \eqref{asymp-s-r-2} for the other components $\bar\gamma^{\pm}_{j}$ of the contour $\Sigma_{R}$. In consequence, we obtain the existence of a constant $c>0$ such that, for any $1\le j\le 4$, we have the following asymptotic behavior
\begin{align}\label{asymp-s-r-3}
\|I - S_{R}\|_{(L^{2}\cap L^{\infty})(\bar\gamma_{j}^{\pm})} = O(e^{-ct}),\quad t\to\infty,
\end{align}
which, in particular, leads to \eqref{bbb12}. Furthermore, combining \eqref{eq-mm-nn-1}, \eqref{ee-m}, \eqref{ee-4-gg}, \eqref{ee-4} and \eqref{asymp-s-r-3} yields \eqref{bbb1} and the proof of the proposition is completed. \hfill $\square$

\begin{lemma}
There is $t_{1}>0$ such that, for any $t>t_{1}$, the RH problem $(R1)-(R3)$ defined on the contour $\Sigma_{R}$ admits a unique solution $R(t,z)$ with the property that 
\begin{align}\label{aaa1}
\|R_{-} - I\|_{L^{2}(\Sigma_{R})} & = O(t^{-1/2}), && \hspace{-40pt} t\to\infty,\\ \label{aaa1b}
\|R_{-} - I - \RR I\|_{L^{2}(\Sigma_{R})} & = O(t^{-1}),&&\hspace{-40pt} t\to\infty.
\end{align}
\end{lemma}
\proof Using Lemma \ref{lem-est-3} we obtain the existence of $t_{0}>0$ such that 
\begin{align}\label{est-t9c}
\|S_{R} - I\|_{(L^{2}\cap L^{\infty})(\Sigma_{R})}\lesssim t^{-1/2}, \quad t\ge t_{0}.
\end{align}
Let us take arbitrary $\rho\in L^{2}_{I}(\Sigma_{R})$, where $\rho = \rho_{0}+\rho_{\infty}$ for $\rho_{0}\in L^{2}(\Sigma_{R})$ and $\rho_{\infty}\in M_{2\times 2}(\C)$. Then, by the linearity of the Cauchy operator, we have
\begin{equation*}
\RR(\rho)= \CC_{-}((\rho-\rho_{\infty})(S_{R} - I))+ \CC_{-}(\rho_{\infty}(S_{R} - I)).
\end{equation*}
Therefore $\RR(\rho)\in L^{2}( \Sigma_{R})$ and the following estimates hold
\begin{equation}
\begin{aligned}\label{est-t5c}
\|\RR(\rho)\|_{L^{2}(\Sigma_{R})}&\lesssim \|\rho_{0}(S_{R} - I)\|_{L^{2}(\Sigma_{R})}+
\|\rho_{\infty}\|\|S_{R} - I\|_{L^{2}(\Sigma_{R})}\\
&\le \|S_{R}-I\|_{(L^{2}\cap L^{\infty})(\Sigma_{R})}\left(\|\rho_{0}\|_{L^{2}(\Sigma_{R})}+\|\rho_{\infty}\|\right)\\
&= \|\rho\|_{L^{2}_{I}(\Sigma_{R})}\|S_{R}-I\|_{(L^{2}\cap L^{\infty})(\Sigma_{R})}.
\end{aligned}
\end{equation}
Combining the inequalities \eqref{est-t9c} and \eqref{est-t5c}, gives
\begin{align}\label{ee10c}
\|\RR(\rho)\|_{L^{2}_{I}(\Sigma_{R})}\lesssim t^{-1/2}\|\rho\|_{L^{2}_{I}(\Sigma_{R})}, \quad t>t_{0},
\end{align}
which, in particular, implies that 
\begin{align}\label{ee11c}
\|\RR\|_{L^{2}_{I}(\Sigma_{R})}\lesssim t^{-1/2}, \quad \|\RR I\|_{L^{2}_{I}(\Sigma_{R})} \lesssim t^{-1/2}, \quad t>t_{0}.
\end{align}
Furthermore \eqref{ee10c} shows that there is $t_{1}>t_{0}$ such that 
\begin{equation*}
\|\RR\|_{L^{2}_{I}(\Sigma_{R})}<1/4,\quad t> t_{1}.
\end{equation*}
Consequently the equation $\rho - \RR(\rho) = I$ has a unique solution $\rho\in L^{2}_{I}( \Sigma_{R})$, given by the Neumann series
$\rho = \sum_{i=0}^{\infty} \RR^{i}I$, which is convergent in the space $L^{2}_{I}(\Sigma_{R})$. Taking into account \eqref{ll-2} and the inequalities \eqref{ee10c}, \eqref{ee11c} yields
\begin{equation*}
\begin{aligned}
&\|R_{-} - I\|_{L^{2}_{I}(\Sigma_{R})}= \|\rho - I\|_{L^{2}_{I}(\Sigma_{R})}\le \sum_{i=1}^{\infty} \|\RR^{i}I\|_{L^{2}_{I}(\Sigma_{R})}\\
&\qquad \le \|\RR I\|_{L^{2}_{I}(\Sigma_{R})}\sum_{i=0}^{\infty} \|\RR\|^{i}_{L^{2}_{I}(\Sigma_{R})}
\le\|\RR I\|_{L^{2}_{I}(\Sigma_{R})} \lesssim t^{-1/2},\quad t>t_{1},
\end{aligned}
\end{equation*}
which proves \eqref{aaa1}. On the other hand, using \eqref{ll-2} and the inequalities \eqref{ee11c}, we obtain the following estimates
\begin{equation*}
\begin{aligned}
&\|R_{-}-I-\RR I\|_{L^{2}_{I}(\Sigma_{R})}=\|\rho - I - \RR I\|_{L^{2}_{I}(\Sigma_{R})}\le \sum_{i=2}^{\infty} \|\RR^{i}I\|_{L^{2}_{I}(\Sigma_{R})}\\
&\quad \le \|\RR\|_{L^{2}_{I}(\Sigma_{R})}\|\RR I\|_{L^{2}_{I}(\Sigma_{R})}\sum_{i=0}^{\infty} \|\RR\|^{i}_{L^{2}_{I}(\Sigma_{R})}
\le\|\RR\|_{L^{2}_{I}(\Sigma_{R})}\|\RR I\|_{L^{2}_{I}(\Sigma_{R})} \lesssim t^{-1}, \ \ t>t_{1},
\end{aligned}
\end{equation*}
that provide \eqref{aaa1b} and complete the proof of the proposition. \hfill $\square$

\section{Proof of Theorem \ref{th-asymptotic-real}}

We begin with the following proposition.
\begin{proposition}\label{prop-asy-1}
We have the following asymptotic relation
\begin{align}\label{kk-ll-22bb}
v(x) = -\frac{\sqrt{-x}}{\pi i}\int_{C}S_R(z')_{12}\,dz' + O((-x)^{-\frac{7}{4}}),\quad x\to-\infty,
\end{align}
where we define $C:=C_{+}\cup C_{-}\cup C_{0}$.
\end{proposition}
\proof Let us observe that using \eqref{u-lim-2}, \eqref{ll-1} and \eqref{ll-2}, we obtain 
\begin{align}\label{rr-kk-1}
v(x) = 2\sqrt{-x}\lim_{z\to\infty}(zR_{12}(z)) = -\frac{\sqrt{-x}}{\pi i}  \int_{\Sigma_R}\left(R_-(z')(S_R(z')-I)\right)_{12}\,dz'.
\end{align}
By the use of \eqref{from-mult} and the H\"older inequality, we have
\begin{align*}
&\left|\int_{\Sigma'_{R}}(R_{-}(z')(S_R(z')-I))_{12}\,dz'\right|\le \int_{\Sigma'_{R}}\|R_{-}(z')(S_R(z')-I)\|\,|dz'|\\
&\qquad\qquad\le \|R_{-}\|_{L^{2}(\Sigma'_{R})} \|S_R-I\|_{L^{2}(\Sigma'_{R})}\le \|R_{-}\|_{L^{2}(\Sigma_{R})} \|S_R-I\|_{L^{2}(\Sigma'_{R})},
\end{align*}
which together with \eqref{bbb12} and \eqref{aaa1} implies that 
\begin{align}\label{rr-kk-2}
\left|\int_{\Sigma'_{R}}(R_{-}(z')(S_R(z')-I))_{12}\,dz'\right| = O(e^{-ct}),\quad t\to\infty,
\end{align}
for some $c>0$. Let us consider the following decomposition
\begin{equation}
\begin{aligned}\label{aa-nn-1}
&\hspace{-5pt}\int_{C}\left(R_-(z')(S_R(z')-I)\right)_{12}\,dz' \!=\! \int_{C}\left((R_{-}(z') \!-\! I \!-\! \RR I(z'))(S_R(z')-I)\right)_{12}\,dz'\\
& \qquad+\int_{C}(S_R(z')-I)_{12}\,dz' +\int_{C}\left(\RR I(z')(S_R(z')-I)\right)_{12}\,dz'.
\end{aligned}
\end{equation}
Using the H\"older inequality and \eqref{from-mult} once again, we obtain
\begin{equation*}
\begin{aligned}
&\left|\int_{C}\left((R_{-}(z') - I - \RR I(z'))(S_R(z')-I)\right)_{12}\,dz'\right|\\
&\ \ \le \int_{C}\|(R_{-}(z') - I - \RR I(z'))(S_R(z')-I)\|\,dz'\\
&\ \ \le \|R_{-} - I - \RR I\|_{L^{2}(C)}\|S_R-I\|_{L^{2}(C)}
\le \|R_{-} - I - \RR I\|_{L^{2}(\Sigma_{R})}\|S_R-I\|_{L^{2}(\Sigma_{R})},
\end{aligned}
\end{equation*}
which combined with \eqref{bbb1} and \eqref{aaa1b} gives
\begin{align*}
\left|\int_{C}\left((R_{-}(z') - I - \RR I(z'))(S_R(z')-I)\right)_{12}\,dz'\right| = O(t^{-3/2}),\quad t\to\infty.
\end{align*}
Therefore, by \eqref{rr-kk-1}, \eqref{rr-kk-2} and \eqref{aa-nn-1}, we infer that 
\begin{align}\label{kk-ll-11}
\hspace{-10pt} v(x) = -\frac{\sqrt{-x}}{\pi i}\left(\int_{C}S_R(z')_{12}+\left(\RR I(z')(S_R(z')-I)\right)_{12}\,dz'\right) 
+ O((-x)^{-\frac{7}{4}}),
\end{align}
as $x\to-\infty$. Let us observe that, by \eqref{asympt-1} and \eqref{asympt-2}, we have
\begin{align}\label{eq-jj-kk}
S_{R}(z) - I = t^{-\frac{1}{2}}F_{\pm}(z) + t^{-1}G_{\pm}(x) + O(t^{-\frac{3}{2}}),\quad t\to\infty.
\end{align}
uniformly for $z\in \partial D(z_{\pm},\delta)$, where $F_{\pm}(z)$ and $G_{\pm}(z)$ are functions given by 
\begin{align*}
F_{+}(z) & := \begin{pmatrix} 0 & \frac{-\nu s_{3}}{h_{1}}e^{\frac{2it}{3}}\frac{\beta(z)^{2}}{\zeta(z)} \\[7pt] 
\frac{-h_{1}}{s_{3}}e^{-\frac{2it}{3}}\frac{\beta(z)^{-2}}{\zeta(z)} & 0 \end{pmatrix}, \ \ G_{+}(z):=\begin{pmatrix} \frac{\nu(\nu+1)}{2\zeta(z)^{2}} & 0 \\[7pt] 0 & -\frac{\nu(\nu-1)}{2\zeta(z)^{2}} \end{pmatrix},
\end{align*}
\begin{align*}
F_{-}(z) & :=\begin{pmatrix}  0 & \frac{h_{1}}{s_{3}}e^{-\frac{2it}{3}}\frac{\beta(-z)^{-2}}{\zeta(-z)} \\[7pt] \frac{\nu s_{3}}{h_{1}}e^{\frac{2it}{3}}\frac{\beta(-z)^{2}}{\zeta(-z)} & 0 \end{pmatrix},  \ \ G_{-}(z) :=\begin{pmatrix}  -\frac{\nu(\nu-1)}{2\zeta(-z)^{2}} & 0\\[7pt] 0 & \frac{\nu(\nu+1)}{2\zeta(-z)^{2}} \end{pmatrix}.
\end{align*}
Let us express the term $\RR I$ in the following form 
\begin{equation}
\begin{aligned}\label{kk-bb}
&\RR I = \CC_{-}(S_{R} - I) = \CC_{-}[(S_{R} - I)\chi_{C_{+}}]+ \CC_{-}[(S_{R} - I)\chi_{C_{-}}]\\
&\qquad + \CC_{-}[(S_{R} - I)\chi_{C_{0}}]+ \CC_{-}[(S_{R} - I)\chi_{\Sigma'_{R}}]\\
& = t^{-\frac{1}{2}}\CC_{-}[F_{+}\chi_{C_{+}}] + t^{-\frac{1}{2}}\CC_{-}[F_{-}\chi_{C_{-}}]+ \CC_{-}[(S_{R} - I)\chi_{C_{0}}]+ \CC_{-}[(S_{R} - I)\chi_{\Sigma'_{R}}]\\
& \qquad +\CC_{-}[(S_{R} - I - t^{-\frac{1}{2}}F_{+})\chi_{C_{+}}]+ \CC_{-}[(S_{R} - I - t^{-\frac{1}{2}}F_{-})\chi_{C_{-}}]
\end{aligned}
\end{equation}
and take arbitrary $j\in\{+,-,0\}$. In view of the inequality \eqref{ineq-c} with $P_{1}=\Sigma'_{R}$ and $P_{2}=C_{j}$, we obtain  
\begin{align*}
\|\CC_{-}[(S_{R} - I)\chi_{\Sigma'_{R}}]\|_{L^{2}(C_{j})}\le \|\CC_{-}\|_{L^{2}(\Sigma_{R})} 
\|S_{R} - I\|_{L^{2}(\Sigma'_{R})},\quad t>0,
\end{align*}
which together with \eqref{bbb12} implies the existence of $c>0$ such that 
\begin{align}\label{kk-bb-1}
\|\CC_{-}[(S_{R} - I)\chi_{\Sigma'_{R}}]\|_{L^{2}(C_{j})} = O(e^{-ct}),\quad t\to\infty.
\end{align}
On the other hand, using \eqref{ineq-c} with $P_{1}=C_{\pm}$ and $P_{2}=C_{j}$, gives
\begin{align*}
\|\CC_{-}[(S_{R} - I - t^{-\frac{1}{2}}F_{\pm})\chi_{C_{\pm}}]\|_{L^{2}(C_{j})}\le \|\CC_{-}\|_{L^{2}(\Sigma_{R})} \|S_{R} - I - t^{-\frac{1}{2}}F_{\pm}\|_{L^{2}(C_{\pm})}
\end{align*}
and therefore, taking into account \eqref{eq-jj-kk}, we obtain the asymptotic relation
\begin{align}\label{kk-bb-2}
\|\CC_{-}(S_{R} - I - t^{-\frac{1}{2}}F_{\pm})\|_{L^{2}(C_{j})}=O(t^{-1}),\quad t\to\infty.
\end{align}
Similarly, applying the inequality \eqref{eq-mm-nn-1} with $P_{1}=C_{0}$ and $P_{2}=C_{j}$, we deduce that
\begin{align*}
\|\CC_{-}[(S_{R} - I)\chi_{C_{0}}]\|_{L^{2}(C_{j})}\le \|\CC_{-}\|_{L^{2}(\Sigma_{R})}\|S_{R} - I\|_{L^{2}(C_{0})}
\end{align*}
which together with \eqref{bbb12kk} provides the relation
\begin{align}\label{kk-bb-3}
\|\CC_{-}[(S_{R} - I)\chi_{C_{0}}]\|_{L^{2}(C_{j})} = O(t^{-1}),\quad t\to\infty.
\end{align}
Then \eqref{kk-bb}, \eqref{kk-bb-1}, \eqref{kk-bb-2} and \eqref{kk-bb-3}, imply
\begin{align}\label{kk-bb-4}
\RR I  = t^{-\frac{1}{2}}\CC_{-}(F_{+}\chi_{C_{+}}) + t^{-\frac{1}{2}}\CC_{-}(F_{-}\chi_{C_{-}}) + O_{L^{2}(C_{j})}(t^{-1}),\quad t\to\infty,
\end{align}
which together with \eqref{eq-jj-kk} gives
\begin{equation}
\begin{aligned}\label{kk-bb-5}
&\int_{C_{\pm}}\left(\RR I(z')(S_R(z')-I)\right)_{12}\,dz' \\
& \ \  = t^{-\frac{1}{2}}\int_{C_{\pm}}\left[(\CC_{-}(F_{+}\chi_{C_{+}})+\CC_{-}(F_{-}\chi_{C_{-}}))(S_R(z')-I)\right]_{12}\,dz'+O(t^{-\frac{3}{2}})\\
& \ \  = t^{-1}\int_{C_{\pm}}\left[(\CC_{-}(F_{+})(z') + \CC_{-}(F_{-})(z'))F_{\pm}(z')\right]_{12}\,dz' + O(t^{-\frac{3}{2}}).
\end{aligned}
\end{equation}
From the definition of the functions $F_{\pm}$, we find that the matrix $\CC_{-}(F_{+}) + \CC_{-}(F_{-})$ has 
the following form 
\begin{align*}
\CC_{-}(F_{+})(z') + \CC_{-}(F_{-})(z') = \begin{pmatrix} 0 & * \\[2pt] * & 0 \end{pmatrix},\quad z'\in \partial D(z_{+},\delta)\cup \partial D(z_{-},\delta).
\end{align*}
This implies that
\begin{align*}
\left((\CC_{-}(F_{+})(z') + \CC_{-}(F_{-})(z'))F_{\pm}(z')\right)_{12} = 0,\quad z'\in \partial D(z_{+},\delta)\cup \partial D(z_{-},\delta).
\end{align*}
and consequently, by \eqref{kk-bb-5}, we have
\begin{align}\label{mm-nn-1}
\int_{C_{\pm}}\left(\RR I(z')(S_R(z')-I)\right)_{12}\,dz' = O(t^{-\frac{3}{2}}),\quad t\to\infty.
\end{align}
On the other hand, the inequality 
\begin{align*}
&\left|\int_{C_{0}}\left(\RR I(z')(S_R(z')-I)\right)_{12}\,dz'\right|\le \int_{C_{0}}\|\RR I(z')(S_R(z')-I)\|\,dz'\\
&\hspace{60pt}\le \|\RR I\|_{L^{2}(C_{0})}\|S_R-I\|_{L^{2}(C_{0})}
\end{align*}
and the asymptotic relations \eqref{bbb12kk}, \eqref{kk-bb-4} provide
\begin{align}\label{mm-nn-2}
\int_{C_{0}}\left(\RR I(z')(S_R(z')-I)\right)_{12}\,dz' = O(t^{-\frac{3}{2}}), \quad t\to\infty.
\end{align}
Combining \eqref{kk-ll-11}, \eqref{mm-nn-1} and \eqref{mm-nn-2} we conclude the asymptotic \eqref{kk-ll-22bb} and the proof of the proposition is completed. \hfill $\square$ \\

In the following proposition we derive the contribution to the asymptotic \eqref{asym-1a} coming from the part of the graph $\Sigma_{T}$ located in the neighborhood of the origin. 

\begin{proposition}\label{prop-c-0}
The following asymptotic relation holds
\begin{align}\label{asym-tt-pp}
\int_{C_{0}}S_R(z')_{12}\,dz' = \pi i\alpha t^{-1} + O(t^{-2}),\quad t\to \infty.
\end{align}
\end{proposition}
\proof In view of the definition of the jump matrix $S_{R}(z)$ and \eqref{equa-p}, we have
\begin{align*}
S_{R}(z)=T^{(0)}(z)N(z)^{-1} = E(z)\bar \Phi(t\eta(z))e^{-it\eta(z)\sigma_3}E(z)^{-1},\quad z\in\Sigma_{R}.
\end{align*}
Observe that, by the point $(d)$  of Theorem \ref{th-l}, the function $\bar \Phi(z)$ has the following asymptotic behavior
\begin{align*}
\bar \Phi(z) = \left(I - \frac{\alpha}{2z}\sigma_1 + H(z)\right)e^{iz\sigma_3}, \quad z\to\infty.
\end{align*}
If $v(x;\alpha,k)$ is a real Ablowitz-Segur solution of the PII equation, then 
\begin{equation}
\begin{aligned}\label{asum-s-r}
S_{R}(z)& = E(z)\left(I - \frac{\alpha}{2t\eta(z)}\sigma_1 + H(t\eta(z))\right)E(z)^{-1}\\
& = I - \frac{\alpha}{2t\eta(z)}\begin{pmatrix} 0 & [E_{11}(z)]^{2} \\[3pt] [E_{11}(z)]^{-2} & 0\end{pmatrix} + 
E(z)H(t\eta(z))E(z)^{-1},
\end{aligned}
\end{equation}
where the matrix coefficient $E_{11}(z)$ is given by $$E_{11}(z) = \left(\frac{z+1/2}{1/2-z}\right)^{\nu}, \quad z\not\in (-\infty,-1/2]\cup[1/2,\infty).$$ In the above formula the branch cut is chosen such that $\mathrm{arg}\,(1/2\pm z)\in (-\pi,\pi)$. By the definition of the map $\eta(z)$ (see \eqref{d-zeta}), there is $c_{0}>0$ such that $|\eta(z)|>c_{0}$ for $|z|=\delta$. Since $H_{lm}(z) = O(z^{-2})$ as $z\to\infty$, for $1\le l,m\le 2$, it follows that 
\begin{align}\label{aa-bb-11}
H(t\eta(z)) = O(t^{-2}),\quad t\to\infty,
\end{align}
uniformly for $z\in\partial D(0,\delta)$. On the other hand, the fact that the function $E(z)$ is holomorphic and invertible in the neighborhood of the origin, implies the existence of a constant $M>0$ such that 
\begin{align*}
\|E(z)\|\le M \quad\text{and}\quad \|E(z)^{-1}\|\le M\quad\text{for} \ \ z\in\partial D(0,\delta),
\end{align*}
which together with \eqref{aa-bb-11} gives 
\begin{align}\label{asum-s-r-2}
E(z)H(t\eta(z))E(z)^{-1} = O(t^{-2}),\quad t\to\infty, 
\end{align}
uniformly for $z\in\partial D(0,\delta)$. Observe that combining \eqref{asum-s-r} and \eqref{asum-s-r-2}, we obtain
\begin{align}\label{eq-11-22}
\int_{C_{0}}S_R(z)_{12}\,dz = - \frac{\alpha}{2t}\int_{C_{0}}\frac{[E_{11}(z)]^{2}}{\eta(z)} \,dz  + O(t^{-2}).
\end{align}
Using the residue method in calculating the integral along the curve $C_{0}$, gives
\begin{equation}
\begin{aligned}\label{eq-11-22-33}
\int_{C_{0}}\frac{[E_{11}(z)]^{2}}{\eta(z)} \,dz &= -2\pi i \,\underset{z=0}{\mathrm{Res}}\,\left(\frac{[E_{11}(z)]^{2}}{\eta(z)}\right) = 
-2\pi i\lim_{z\to 0} \frac{[E_{11}(z)]^{2}}{\eta(z)/z}\\
& = -2\pi i\lim_{z\to 0} \frac{[E_{11}(z)]^{2}}{\eta(z)/z}=-2\pi i\lim_{z\to 0} \frac{[E_{11}(z)]^{2}}{1-4z^{2}/3} = -2\pi i.
\end{aligned}
\end{equation}
Therefore, by \eqref{eq-11-22} and \eqref{eq-11-22-33} we deduce the asymptotic relation \eqref{asym-tt-pp} and the proof of the proposition is completed. \hfill $\square$ \\

In the following two propositions we calculate the contribution to the asymptotic \eqref{asym-1a} coming from the part of the graph $\Sigma_{T}$ located in the neighborhoods of the stationary points $z_{\pm}=\pm1/2$. 

\begin{proposition}\label{prop-asym-rel-1}
We have the following asymptotic relation as $t\to\infty$:
\begin{align}\label{asym-pp-11}
\int_{C_{+}\cup C_{-}}\!\!S_R(z')_{12}\,dz' = -i\pi dt^{-\frac{1}{2}}\cos\left(\frac{2}{3}t - \frac{3}{4}d^{2}\ln(t^{2/3})+\phi\right)+O(t^{-\frac{3}{2}})
\end{align}
where the constants $d$ and $\phi$ are given by the connection formulas \eqref{conn-f-real-1} and \eqref{conn-f-real-2}.
\end{proposition}
\proof By the asymptotic relation \eqref{eq-jj-kk}, the following holds as $t\to\infty$:
\begin{equation}
\begin{aligned}\label{eq-jj-kk-11}
&\hspace{-5pt}\int_{C_{+}\cup C_{-}}S_R(z')_{12}\,dz' \!=\! t^{-\frac{1}{2}}\!\!\int_{C_{+}}\!\!F_{+}(z')_{12}\,dz'+t^{-\frac{1}{2}}\!\!\int_{C_{-}}\!\!F_{-}(z')_{12}\,dz'+O(t^{-\frac{3}{2}})\\
&\hspace{-5pt}\quad = -t^{-\frac{1}{2}}\frac{\nu s_{3}}{h_{1}}e^{\frac{2it}{3}}\!\int_{C_{+}}\! \frac{\beta(z')^{2}}{\zeta(z')}\,dz' + t^{-\frac{1}{2}}\frac{h_{1}}{s_{3}}e^{-\frac{2it}{3}}\!\int_{C_{-}} \!\frac{\beta(-z')^{-2}}{\zeta(-z')} \,dz' + O(t^{-\frac{3}{2}}).
\end{aligned}
\end{equation}
If $v(x;\alpha,k)$ is a real Ablowitz-Segur solution, then the numbers $s_{1},s_{3}$ defined in \eqref{stokes-11bb}, \eqref{stokes-11bbc} are such that $s_{1} = \overline{s_{3}}$. Therefore the results of \cite[Page 28]{MR3670014} say that
\begin{equation}
\begin{aligned}\label{eq-jj-kk-22}
&-t^{-\frac{1}{2}}\frac{\nu s_{3}}{h_{1}}e^{\frac{2it}{3}}\int_{C_{+}} \frac{\beta(z')^{2}}{\zeta(z')}\,dz' + t^{-\frac{1}{2}}\frac{h_{1}}{s_{3}}e^{-\frac{2it}{3}}\int_{C_{-}} \frac{\beta(-z')^{-2}}{\zeta(-z')} \,dz'\\
&\qquad = -i\pi dt^{-\frac{1}{2}}\cos(\frac{2}{3}t - \frac{3}{4}d^{2}\ln(t^{2/3})+\phi),\qquad t>0,
\end{aligned}
\end{equation}
where the constants $d$ and $\phi$ are given by the formulas \eqref{conn-f-real-1} and \eqref{conn-f-real-2}. Consequently, by \eqref{eq-jj-kk-11} and \eqref{eq-jj-kk-22}, we obtain the relation \eqref{asym-pp-11} and the proof of the proposition is completed. \hfill $\square$\\

\noindent{\em Proof of Theorem \ref{th-asymptotic-real}.} If $v(x;\alpha,k)$ is a real Ablowitz-Segur solution of the inhomogeneous PII equation, then 
applying Propositions \ref{prop-c-0} and \ref{prop-asym-rel-1} we obtain 
\begin{align*}
\int_{C}S_R(z')_{12}\,dz' = \pi i\alpha t^{-1} -i\pi dt^{-\frac{1}{2}}\cos\left(\frac{2}{3}t - \frac{3}{4}d^{2}\ln(t^{2/3})+\phi\right) \!+ \!O(t^{-\frac{3}{2}}), \ t\to\infty,
\end{align*}
where $d$ and $\phi$ are given by \eqref{conn-f-real-1}, \eqref{conn-f-real-2}. Substituting this asymptotic relation into \eqref{kk-ll-22bb} and using \eqref{change-var} we obtain the relation \eqref{asym-1a} and the proof of the theorem is completed. \hfill $\square$

\section{Proof of Theorem \ref{th-kdv}}
We begin with the following proposition.
\begin{proposition}\label{prop-lim-d-p}
Assume that $v(x;\alpha,k)$ is real Ablowitz-Segur solution of the PII equation. Then the following limit holds
\begin{align}\label{conv-1a}
\lim_{t\to 0^{+}}t^{-\frac{1}{3}} v(xt^{-\frac{1}{3}};\alpha,k) = c(\alpha,k)\,\delta(x)+\alpha\,\mathrm{p.v.}\,(1/x)
\end{align}
in the space $\SS'(\R)$, where we define 
\begin{equation}\label{t-int}
c(\alpha,k):=\lim_{y\to\infty}\int_{-y}^{y}v(x;\alpha,k)\,dx.
\end{equation}
Furthermore, we have the pointwise limit 
\begin{align}\label{eq-lim-33}
\lim_{\xi\to 0^{\pm}}\widehat{v}(\xi;\alpha,k) = c(\alpha,k)\mp i\pi\alpha.
\end{align}
\end{proposition}
\proof Let us define $\Psi(x):=\frac{2}{3}(-x)^{\frac{3}{2}}-\frac{3}{4}d^{2}\ln(-x)+\phi$, where the constants $d$ and $\phi$ are given by the connection formulas \eqref{conn-f-real-1} and \eqref{conn-f-real-2}. Assume that $x_{0}>1$ is sufficiently large such that
\begin{equation}\label{non-zero}
\Psi'(x)=x^{1/2}-\frac{3}{4}d^2 x^{-1}\neq 0,\quad x\ge x_{0}
\end{equation}
and define the following auxiliary functions
\begin{equation*}
\begin{gathered}
f(x):= \alpha x^{-1}\chi_{\mathbb{R}\setminus(-1,\,1)}(x), \ \
g(x):= d\chi_{(-\infty,\,-x_{0})}(x)(-x)^{-1/4}\cos\Psi(-x), \\ h(x):=v(x)-f(x)-g(x),\quad x\in\R.
\end{gathered}
\end{equation*}
Then, for any $\varphi\in\mathcal{S}(\R)$, $t>0$ and $z\ge x_{0}$, we have the following decomposition
\begin{equation*}
\begin{aligned}
&\int_\R t^{-\frac{1}{3}} v(xt^{-\frac{1}{3}})\varphi(x)\,dx = \int_{zt^{\frac{1}{3}}}^{\infty} t^{-\frac{1}{3}} f(xt^{-\frac{1}{3}})\varphi(x)\,dx+ \int_{zt^{\frac{1}{3}}}^{\infty} t^{-\frac{1}{3}} h(xt^{-\frac{1}{3}})\varphi(x)\,dx \\
&\qquad + \int^{zt^{\frac{1}{3}}}_{-zt^{\frac{1}{3}}} t^{-\frac{1}{3}}v(xt^{-\frac{1}{3}})\varphi(x)\,dx + \int_{-\infty}^{-zt^{\frac{1}{3}}} t^{-\frac{1}{3}}f(xt^{-\frac{1}{3}})\varphi(x)\,dx\\
&\qquad\qquad +\int_{-\infty}^{-zt^{\frac{1}{3}}} t^{-\frac{1}{3}}g(xt^{-\frac{1}{3}})\varphi(x)\,dx + \int_{-\infty}^{-zt^{\frac{1}{3}}} t^{-\frac{1}{3}}h(xt^{-\frac{1}{3}})\varphi(x)\,dx.
\end{aligned}
\end{equation*}
Changing the variables and rearranging the terms of the above sum, we obtain
\begin{equation}
\begin{aligned}\label{eq-11}
&\int_\R t^{-\frac{1}{3}} v(xt^{-\frac{1}{3}})\varphi(x)\,dx = \int_{z}^{\infty} h(x)\varphi(xt^{1/3})\,dx + \int^{-z}_{-\infty} h(x)\varphi(xt^{1/3})\,dx \\
&\quad +\int_{-z}^{z} v(x) \varphi(xt^{1/3})\,dx +\alpha\int_{|x|>zt^{1/3}} \frac{\varphi(x)}{x}\,dx + \int_{-\infty}^{-z} g(x)\varphi(xt^{1/3})\,dx.
\end{aligned}
\end{equation}
It is not difficult to check that
\begin{align}\label{gg-3}
\lim_{t\to 0^{+}}\int_{-z}^{z} v(x)\varphi(xt^{1/3}) \,dx = \varphi(0)\int_{-z}^{z} v(x)\,dx.
\end{align}
Let us observe that the relation \eqref{asym-1a} gives the following asymptotics
\begin{equation}
\begin{aligned}\label{asym-11a}
&h(x)\varphi(xt^{1/3}) = O(x^{-4}),\quad x\to\infty,\\
&h(x)\varphi(xt^{1/3}) = O((-x)^{-7/4}),\quad  x\to-\infty,
\end{aligned}
\end{equation}
that are uniform with respect to $t\in[0,1]$. Hence, by \eqref{asym-11a} and dominated convergence theorem, we obtain
\begin{align}\label{gg-1}
\lim_{t\to 0^{+}}\int_{z}^{\infty}h(x)\varphi(xt^{1/3})\,dx=\varphi(0)\int_{z}^{\infty}h(x)\,dx = O(z^{-3}),\quad z\to\infty
\end{align}
and furthermore
\begin{align}\label{gg-2}
\lim_{t\to 0^{+}}\int^{-z}_{-\infty} h(x)\varphi(xt^{1/3})\,dx = \varphi(0)\int^{-z}_{-\infty} h(x)\,dx = O(z^{-3/4}),\quad z\to\infty.
\end{align}
We proceed to the estimates for the last term of \eqref{eq-11}. To this end, we integrate by parts and observe that, for any $z\ge x_{0}$, we have
\begin{equation}
\begin{aligned}\label{eq-33}
\int_{-\infty}^{-z} g(x)\varphi(xt^{1/3})\,dx &= -\frac{d\sin\Psi(z) \varphi(-zt^{1/3})}{z^{1/4}\Psi'(z)}
+d \int_{z}^{\infty} J_{2}(x)\varphi(-xt^{1/3})\,dx  \\
& \quad + dt^{1/3}\int_{z}^{\infty}J_{1}(x)\varphi'(-xt^{1/3})\,dx,
\end{aligned}
\end{equation}
where we used the notation
\begin{align*}
J_{1}(x):=\frac{\sin\Psi(x)}{x^{1/4}\Psi'(x)} \quad \text{ and } \quad J_{2}(x):=\frac{(x^{1/4}\Psi'(x))'}{x^{1/2}\Psi'(x)^{2}}\sin\Psi(x),\quad x\ge x_{0}.
\end{align*}
In view of \eqref{non-zero}, it is not difficult to check the asymptotics
\begin{equation}\label{asym-1bc}
x^{1/2}\Psi'(x)^{2}\sim x^{3/2}\quad\text{and}\quad [x^{1/4}\Psi'(x)]'\sim x^{-1/4},\quad x\to\infty,
\end{equation}
that yield the following relation $$J_{2}(x)\varphi(-xt^{1/3}) = O(x^{-7/4}),\quad x\to\infty,$$
uniformly for $t\in[0,1]$. Therefore, by the dominated convergence theorem, we have
\begin{align}\label{kk-11}
\lim_{t\to 0^{+}} \int_{z}^{\infty} J_{2}(x)\varphi(-xt^{\frac{1}{3}})\,dx = \varphi(0)\int_{z}^{\infty}J_{2}(x)\,dx =
O(z^{-\frac{3}{4}}), \ \ z\to \infty.
\end{align}
On the other hand, integrating by parts once again, for any $z\ge x_{0}$, we have
\begin{equation}
\begin{aligned}\label{eq-22}
\int_{z}^{\infty}J_{1}(x)\varphi'(-xt^{1/3})\,dx &=\frac{\cos\Psi(z) \varphi'(-zt^{1/3})}{z^{1/4} \Psi'(z)^2} + J_{4}(x)\varphi'(-xt^{1/3})\,dx\\
&\qquad - \int_{z}^{\infty} t^{1/3}J_{3}(x)\varphi''(-xt^{1/3})\,dx,
\end{aligned}
\end{equation}
where we define the functions
\begin{align*}
J_{3}(x) := \frac{\cos\Psi(x)}{x^{1/4}\Psi'(x)^2}  \quad \text{ and } \quad J_{4}(x) := \frac{[x^{1/4}\Psi'(x)^2]'}{x^{1/2}\Psi'(x)^4}\cos\Psi(x),\quad x\ge x_{0}.
\end{align*}
Using \eqref{non-zero} once again, it is not difficult to check that
\begin{align}\label{asym-2b}
x^{1/4}\Psi'(x)^2\sim x^{5/4},\quad [x^{1/4}\Psi'(x)^2]'\sim x^{1/4} \quad\text{and}\quad x^{1/2}\Psi'(x)^4\sim x^{5/2},
\end{align}
as $x\to\infty$, which gives the following asymptotic relations
\begin{align*}
J_{3}(x)\varphi''(-xt^{1/3}) = O(x^{-5/4})\quad\text{and}\quad J_{4}(x)\varphi'(-xt^{1/3}) = O(x^{-9/4}),\quad x\to\infty,
\end{align*}
uniformly for $t\in [0,1]$. Hence, by \eqref{eq-22} and the dominated convergence theorem
\begin{align}\label{kk-22}
\lim_{t\to 0^{+}}t^{1/3}\int_{z}^{\infty}J_{1}(x)\varphi'(-xt^{1/3})\,dx = 0,\quad z\ge x_{0}.
\end{align}
Since $z^{1/4}\Psi'(z) \sim z^{3/4}$ as $z\to\infty$, it follows that
\begin{align*}
-\frac{d\sin\Psi(z) \varphi(0)}{z^{1/4}\Psi'(z)} = O(z^{-3/4}),\quad z\to\infty,
\end{align*}
which combined with \eqref{eq-33}, \eqref{kk-11} and \eqref{kk-22}, yields
\begin{align}\label{gg-4}
\lim_{t\to 0^{+}}\int_{-\infty}^{-z} g(x)\varphi(xt^{1/3})\,dx = O(z^{-3/4}),\quad z\to\infty.
\end{align}
Passing in the equality \eqref{eq-11} to the limit with $t\to 0^{+}$, for any $z\ge x_{0}$, we obtain
\begin{align*}
&\!\!\!\!\lim_{t\to 0^{+}}\int_\R t^{-\frac{1}{3}} v(xt^{-\frac{1}{3}})\varphi(x)\,dx \!=\! \lim_{t\to 0^{+}}\!\int_{z}^{\infty}\! h(x)\varphi(xt^{\frac{1}{3}})\,dx \!+\! \lim_{t\to 0^{+}}\!\int^{-z}_{-\infty}\! h(x)\varphi(xt^{\frac{1}{3}})\,dx \\
&\qquad +\varphi(0) \int_{-z}^{z} v(x) \,dx +\alpha\lim_{t\to 0^{+}}\int_{|x|>t} \frac{\varphi(x)}{x}\,dx + \lim_{t\to 0^{+}}\int_{-\infty}^{-z} g(x)\varphi(xt^{\frac{1}{3}})\,dx.
\end{align*}
Letting $z\to\infty$ and applying \eqref{gg-3}, \eqref{gg-1}, \eqref{gg-2} and \eqref{gg-4} provide
\begin{align*}
\lim_{t\to 0^{+}}\int_\R t^{-\frac{1}{3}} v(xt^{-\frac{1}{3}})\varphi(x)\,dx=\varphi(0)\lim_{z\to\infty}\int_{-z}^{z} v(x)\,dx+\alpha\lim_{t\to 0^{+}}\int_{|x|>t}\frac{\varphi(x)}{x}\,dx,
\end{align*}
which together with \eqref{t-int} completes the proof of the equality \eqref{conv-1a}. To show that the limit \eqref{eq-lim-33} is also satisfied, we observe that
\begin{align*}
\widehat{v}(0) = \lim_{y\to\infty} \int_{-y}^y v(x)\,dx = c(\alpha,k)
\end{align*}
and we consider the following decomposition
\begin{equation}\label{dec-1}
\widehat{v}(\xi) - c(\alpha,k) = \widehat{v}(\xi)-\widehat{v}(0)=(\widehat{f}(\xi)-\widehat{f}(0))+(\widehat{g}(\xi)-\widehat{g}(0))+(\widehat{h}(\xi)-\widehat{h}(0)).
\end{equation}
By the asymptotic relation \eqref{asym-1a}, the function $h$ is an element of the space $L^{1}(\R)$. Hence its Fourier transform is continuous and, in particular, we have
\begin{equation}\label{s1}
\lim_{\xi\to0}\widehat{h}(\xi)=\widehat{h}(0).
\end{equation}
Since $f$ is an odd function, it is not difficult to check that
\begin{align*}
\widehat{f}(0)= \lim_{y\to\infty}\int_{-y}^{y} f(x)\,dx = 0
\end{align*}
and consequently we obtain
\begin{equation*}
\widehat{f}(\xi)-\widehat{f}(0) = -2i\alpha\lim_{y\to\infty}\int_{1}^{y}\sin(\xi x)x^{-1}\,dx=-2i\alpha\left(\text{sgn}(\xi)\frac{\pi}{2}-\int_{0}^{\xi}\frac{\sin x}{x}\,dx\right).
\end{equation*}
This in turn implies that
\begin{align}\label{s2plus}
\lim_{\xi\to0^{\pm}}(\widehat{f}(\xi)-\widehat{f}(0)) =\mp i\pi\alpha.
\end{align}
In view of the definition of the function $g$, we have
\begin{align*}
\widehat{g}(\xi) - \widehat{g}(0)=d\lim_{y\to\infty}\int_{x_{0}}^{y} (e^{i\xi x}-1)x^{-1/4}\cos\Psi(x)\,dx.
\end{align*}
Therefore, using \eqref{non-zero} and integrating by parts, we obtain
\begin{equation}
\begin{aligned}\label{est-11}
&\widehat{g}(\xi) - \widehat{g}(0)=d\lim_{y\to\infty}\int_{x_{0}}^{y} \frac{e^{i\xi x}-1}{x^{1/4}\Psi'(x)}(\sin\Psi(x))'\,dx \\
& \,= -d\frac{e^{i\xi x_{0}}-1}{c^{1/4}\Psi'(x_{0})}\sin\Psi(x_{0}) -d\lim_{y\to\infty}\int_{x_{0}}^{y}(i\xi e^{i\xi x}J_{1}(x) - (e^{i\xi x}-1)J_{2}(x))\,dx.
\end{aligned}
\end{equation}
On the other hand, using again the asymptotic relation \eqref{asym-1bc}, we have
$$(e^{i\xi x}-1)J_{2}(x) = O(x^{-7/4}),\quad x\to\infty$$ uniformly for $\xi\in[-1,1]$, which by the dominated convergence theorem, gives
\begin{align}\label{eq-i-1}
\lim_{\xi\to 0}\lim_{y\to\infty}\int_{x_{0}}^{y}(e^{i\xi x}-1)J_{2}(x)\,dx = \int_{x_{0}}^{\infty}\lim_{\xi\to 0}(e^{i\xi x}-1)J_{2}(x)\,dx = 0.
\end{align}
Let us observe that using \eqref{non-zero} and integrating by parts once again, we find that
\begin{equation}
\begin{aligned}\label{est-22}
&\lim_{y\to\infty}\int_{x_{0}}^{y}i\xi e^{i\xi x}J_{1}(x)\,dx = \lim_{y\to\infty}\int_{x_{0}}^{y}\frac{i\xi e^{i\xi x}}{x^{1/4}\Psi'(x)}\sin\Psi(x)\,dx \\
&\qquad= \frac{i\xi e^{i\xi x_{0}}}{x_{0}^{1/4}\Psi'(x_{0})^{2}}\cos\Psi(x_{0}) - \lim_{y\to\infty}\int_{x_{0}}^{y}(\xi^{2}e^{i\xi x}J_{3}(x)+i\xi e^{i\xi x}J_{4}(x))\,dx,
\end{aligned}
\end{equation}
Applying the relation \eqref{asym-2b} provides
$$\xi^{2}e^{i\xi x}J_{3}(x) = O(x^{-5/4}) \quad\text{and}\quad i\xi e^{i\xi x}J_{4}(x) = O(x^{-9/4}),\quad  x\to\infty,$$ uniformly for $\xi\in[-1,1]$. This, by the dominated convergence theorem, gives
\begin{align*}
\lim_{\xi\to 0}\lim_{y\to\infty}\int_{x_{0}}^{y}(\xi^{2}e^{i\xi x}J_{3}(x)+i\xi e^{i\xi x}J_{4}(x))\,dx = 0,
\end{align*}
which together with \eqref{est-22} yields
\begin{align}\label{eq-i-2}
\lim_{\xi\to 0}\lim_{y\to\infty}\int_{x_{0}}^{y}i\xi e^{i\xi x}J_{1}(x)\,dx = 0.
\end{align}
Therefore, taking into account \eqref{est-11}, \eqref{eq-i-1} and \eqref{eq-i-2}, we infer that
\begin{equation}\label{eq-cc-11}
\lim_{\xi\to0}(\widehat{g}(\xi)-\widehat{g}(0))=0.
\end{equation}
Combining \eqref{dec-1}, \eqref{s1}, \eqref{s2plus} and \eqref{eq-cc-11} yields the limit \eqref{eq-lim-33} and the proof of the proposition is completed. \hfill $\square$ \\

In the proof of Theorem \ref{th-kdv} we will also need the following result, which provides a formula expressing the value of the Cauchy principal value integral \eqref{t-int} for the real Ablowitz-Segur solutions in the terms of the parameters $\alpha$ and $k$. 
\begin{theorem}\label{th-total}
If $v(x;\alpha,k)$ is a real Ablowitz-Segur solution for the inhomogeneous second Painlev\'e equation, then 
\begin{align}\label{wz1}
&\lim_{x\to\infty} \int_{-x}^x v(y;\alpha,k)\,dy = \frac{1}{2}\ln\left(\frac{\cos(\pi\alpha)+k}{\cos(\pi\alpha)-k}\right).
\end{align}
\end{theorem}
In the case of the homogeneous PII equation $\alpha=0$ and $k\in(-1,1)$, the above formula was established in \cite[Theorem 2.1]{MR2501035} by the application of the steepest descent analysis to the RH problem associated with the PII equation. Similar techniques were used in \cite[Theorem 1.1]{kok} to prove the formula \eqref{wz1} for $\alpha$ and $k$ satisfying the general condition \eqref{stokes-11bbc}.\\[5pt]
\noindent\text{\em Proof of Theorem \ref{th-kdv}.} Given $a\in\R$ and $b\in(-1,1)$, let us write $\alpha:=-b/2$ and take $k\in (-\cos(\pi\alpha),\cos(\pi\alpha))$ such that
\begin{align*}
\frac{1}{2}\ln\left(\frac{\cos(\pi\alpha)+k}{\cos(\pi\alpha)-k}\right) = -\frac{a}{2}.
\end{align*}
By Theorem \ref{th-total}, the real Ablowitz-Segur solution $v=v(\,\cdot\,;\alpha,k)$ of the PII equation \eqref{PII} satisfies the equality
\begin{align}\label{eq-arg}
c(\alpha,k) = \lim_{y\to\infty} \int_{-y}^y v(x;\alpha,k)\,dx = -\frac{a}{2}.
\end{align}
Given $\varphi \in\SS(\R)$, let us consider the scaled function $\tilde\varphi(x):=\varphi(x3^{1/3})$ for $x\in\R$. Then, by the form \eqref{self-sim} of the function $u$ and Proposition \ref{prop-lim-d-p}, we infer that
\begin{align*}
&\lim_{t\to 0^{+}}\int_{\R}u(t,x)\varphi(x)\,dx = -2\lim_{t\to 0^{+}} t^{-1/3}\int_{\R}v(x t^{-1/3})\tilde\varphi(x)\,dx \\
&\ \ = -2c(\alpha,k)\,\tilde\varphi(0)-2\alpha\,\lim_{\ve\to 0^{+}}\int_{|x|>\ve} \frac{\tilde\varphi(x)}{x}\,dx
 =  a\,\varphi(0)+b\lim_{\ve\to 0^{+}}\int_{|x|>\ve} \frac{\varphi(x)}{x}\,dx,
\end{align*}
where the last equality follows from \eqref{eq-arg} and fact that the Dirac delta function and Cauchy principal value are invariant under dilation. This completes the proof of the convergence \eqref{conv-1}. To show the pointwise limit \eqref{eq-lim-2}, we use Proposition \ref{prop-lim-d-p} once again and observe that the limit \eqref{eq-lim-33} gives
\begin{align*}
\lim_{t\to 0^{+}}\widehat{u}(t,\xi) &= -2\lim_{t\to 0^{+}}\widehat{v}(\xi(3t)^{1/3}) = -2c(\alpha,k) +2i\pi\alpha\,\mathrm{sgn}\,\xi 
=a - i\pi b\,\mathrm{sgn}\,\xi,
\end{align*}
and the proof of the theorem is completed. \hfill $\square$ 

\section{Appendix}
\begin{figure}[h]
\begin{subfigure}{0.49\textwidth}
\begin{center}
\includegraphics[scale=0.7]{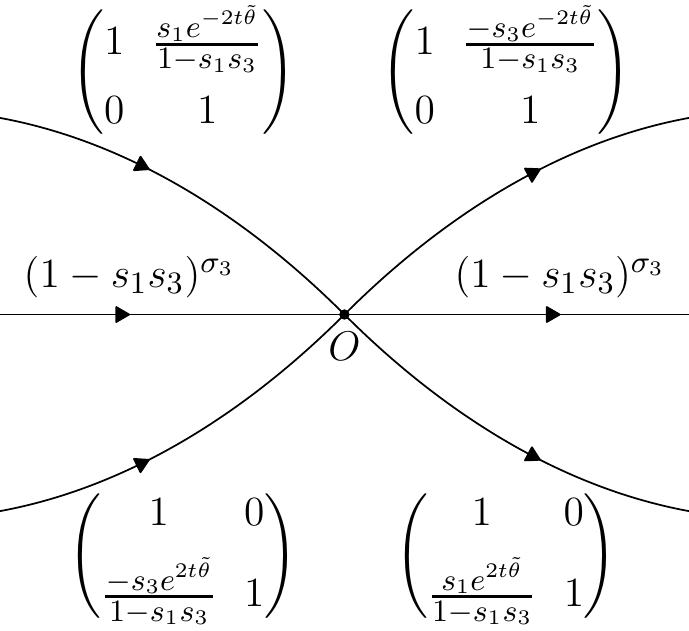}
\end{center}
\end{subfigure}
\begin{subfigure}{0.49\textwidth}
\begin{center}
\includegraphics[scale=0.8]{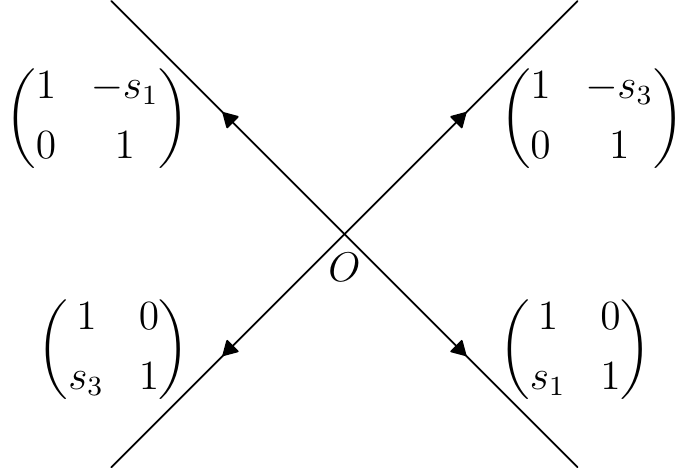}
\end{center}  
\end{subfigure}
\caption{The graph $D(0,\delta)\cap \Sigma_{T}$ for the RH problem fulfilled by the parametrix around the origin.}
\label{fig14}
\end{figure}

In this section we consider local parametrices that are necessary in the analysis of the RH problem (T1)--(T5), considered in Section 2. Let us recall that the triple $(s_{1},s_{2},s_{3})$ and the parameters $\alpha$, $k$ are defined by the formulas \eqref{stokes-11bb} and \eqref{stokes-11bbc}. At the beginning we consider a parametrix near the origin that were obtained in \cite[Section 11.6]{MR2264522}. For this purpose, we define $\bar\Sigma$ to be the contour consisting of the rays $\mathrm{arg}\,z = \pm\pi/4$ and $\mathrm{arg}\,z = \pm 3\pi/4$ (see the right diagram of Figure \ref{fig14}). 
\begin{theorem}\label{th-l}
The function $\bar\Phi(z)$ satisfies the following RH problem.\\[2pt]
\noindent\makebox[5.5mm][l]{$(a)$}\parbox[t][][t]{121mm}{The function $\bar\Phi(z)$ is an analytic function on $\C\setminus\bar\Sigma$;}\\[2pt]
\noindent\makebox[5.5mm][l]{$(b)$}\parbox[t][][t]{121mm}{On the contour $\bar\Sigma$, the following jump relation is satisfied
\begin{equation*}
\bar\Phi_+(z) = \bar\Phi_-(z)\bar S(z),\quad z\in \bar\Sigma,
\end{equation*}
where the jump matrix $\bar S(z)$ is given on the right diagram of Figure \ref{fig14}.}\\[2pt]
\noindent\makebox[5.5mm][l]{$(c)$}\parbox[t][][t]{122mm}{We have the following asymptotic relation 
\begin{align*}
\bar\Phi(z)= O\begin{pmatrix}|z|^{-|\alpha|}& |z|^{-|\alpha|}\\[5pt] |z|^{-|\alpha|}& |z|^{-|\alpha|}\end{pmatrix},\quad z\to 0.
\end{align*}}\\[5pt]
\noindent\makebox[6mm][l]{$(d)$}\parbox[t][][t]{121mm}{The function $\bar\Phi(z)$ has the following behavior at infinity
\begin{equation*}
\bar\Phi(z) =(I - \frac{\alpha}{2z}\sigma_1 + O(z^{-2}))e^{iz\sigma_3}, \quad z\to\infty.
\end{equation*}}
\end{theorem}
Let us assume that the parameter $\nu$ has the following value
\begin{align*}
\nu := -(2\pi i)^{-1}\ln(1-s_1s_3) = -(2\pi i)^{-1}\ln(\cos^2(\pi\alpha)-k^{2})
\end{align*}
and consider the function $N(z)$ given by the formula
\begin{equation}\label{f-z}
N(z):=\left(\frac{z+1/2}{z-1/2}\right)^{\nu\sigma_{3}},\quad z\in \C\setminus[z_{-},z_{+}],
\end{equation}
where $[z_{-},z_{+}]$ is the segment connecting the stationary points $z_{\pm}=\pm \frac{1}{2}$ and the branch cut is taken such that $\mathrm{arg}\,(z\pm1/2)\in(-\pi,\pi)$. Let us assume that $T^{(0)}(z)$ is a function on the ball $D(0,\delta)$, which is given by the formula
\begin{equation}\label{equa-p} 
T^{(0)}(z):=\left\{\begin{aligned}
&E(z)\bar \Phi(t\eta(z))e^{-it\eta(z)\sigma_3} e^{-i\pi\nu\sigma_3},&&\mathrm{Im}\,z>0,\\
&E(z)\bar \Phi(t\eta(z))e^{-it\eta(z)\sigma_3} e^{i\pi\nu\sigma_3},&&\mathrm{Im}\,z<0,
\end{aligned}\right.
\end{equation}
where the function $\eta(z)$ is defined by the formula \eqref{d-eta} and $E(z)$ is a map given by
\begin{equation*}
E(z):= N(z)e^{i\pi\nu\sigma_3},\quad \mathrm{Im}\, z>0\quad\text{and}\quad E(z):=N(z)e^{-i\pi\nu\sigma_3},\quad \mathrm{Im}\, z<0.
\end{equation*}
\begin{theorem}\label{lok-aprox}
The function $T^{(0)}(z)$ is a solution of the following RH problem.\\[3pt]
\noindent\makebox[5.5mm][l]{$(a)$}\parbox[t][][t]{121mm}{The function $T^{(0)}(z)$ is analytic in $D(0,\delta)\setminus\Sigma_T$.}\\[3pt]
\noindent\makebox[5.5mm][l]{$(b)$}\parbox[t][][t]{121mm}{On the contour $D(0,\delta)\cap \Sigma_{T}$ the function $T^{(0)}(z)$ satisfies the same jump conditions as $T(z)$ (see Figure \ref{fig14}).}\\[3pt]
\noindent\makebox[5.5mm][l]{$(c)$}\parbox[t][][t]{121mm}{The function $T^{(0)}(z)$ has the following asymptotic behavior 
\begin{equation*}
T^{(0)}(z)N(z)^{-1} = I + O(t^{-1}),\quad  t\to\infty,
\end{equation*}
uniformly for $z\in \partial D(0,\delta)$.}\\[3pt]
\noindent\makebox[5.5mm][l]{$(d)$}\parbox[t][][t]{122mm}{We have the following asymptotic relation 
\begin{align*}
T^{(0)}(z)= O\begin{pmatrix}|z|^{-|\alpha|}& |z|^{-|\alpha|}\\[5pt] |z|^{-|\alpha|}& |z|^{-|\alpha|}\end{pmatrix},\quad z\to 0.
\end{align*}}
\end{theorem}

We proceed to the local parametrices around stationary points $z_{\pm}=\pm1/2$ for the deformed RH problem associated with the PII equation. To this end we use the well-known parabolic cylinder functions (see e.g. \cite{ba-er} and \cite{MR0232968} for more details) to define the following matrix-valued holomorphic map  
\begin{align*}
Z_{0}(z):=2^{-\frac{\sigma_{3}}{2}}\begin{pmatrix} D_{-\nu-1}(iz) & D_{\nu}(z) \\[5pt]
\frac{d}{dz} D_{-\nu-1}(iz) & \frac{d}{dz} D_{\nu}(z) \end{pmatrix}\begin{pmatrix} e^{\frac{i\pi}{2}(\nu+1)} & 0 \\[5pt] 0 & 1 \end{pmatrix},
\end{align*}
and define the triangular matrices 
\begin{gather*}
H_{0} \!=\! \begin{pmatrix} 1&  0 \\[2pt] h_{0} & 1\end{pmatrix}, \ \
H_{1} \!=\! \begin{pmatrix} 1&  h_{1} \\[2pt] 0 &  1\end{pmatrix}, \ \
H_{2} \!=\! \begin{pmatrix} 1& \!0 \\[2pt] -h_{0}e^{-2i\pi\nu} & \! 1\end{pmatrix}, \ \
H_{3} \!=\! \begin{pmatrix} 1& \!-h_{1}e^{2i\pi\nu} \\[2pt] 0 &\! 1\end{pmatrix}, 
\end{gather*}
where the complex constants $h_{0}$ and $h_{1}$ are given by 
\begin{align*}
h_{0}:=-i\sqrt{2\pi}\,\Gamma(\nu+1)^{-1},\quad h_{1}:=\sqrt{2\pi}e^{i\pi\nu}\,\Gamma(-\nu)^{-1}.
\end{align*}
Let us assume that $Z(z)$ is a sectionally holomorphic matrix function given by
\begin{equation*}
Z(z):=\left\{\begin{aligned}
&Z_{0}(z), \quad \mathrm{arg}\,z\in(-\tfrac{\pi}{4},0),\\
&Z_{j}(z), \quad \mathrm{arg}\,z\in(\tfrac{(j-1)\pi}{2},\tfrac{j\pi}{2}), \ j=1,2,3,\\
&Z_{4}(z), \quad \mathrm{arg}\,z\in(\tfrac{3\pi}{2},\tfrac{7\pi}{4}),
\end{aligned}\right.
\end{equation*}
where the functions $Z_{j}(z)$, for $1\le j\le 4$, are given by the recurrence relation
\begin{align*}
Z_{j+1}(z) = Z_{j}(z)H_{j},\quad j=0,1,2,3.
\end{align*}

\begin{figure}[h]
\begin{subfigure}{0.49\textwidth}
\begin{center}
\hspace{-10pt}\includegraphics[scale=0.72]{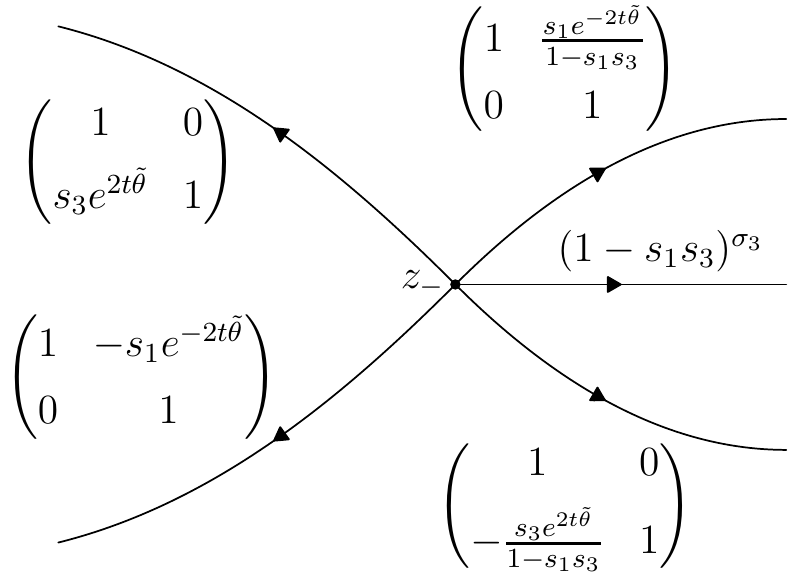}
\end{center}  
\end{subfigure}
\begin{subfigure}{0.49\textwidth}
\begin{center}
\includegraphics[scale=0.72]{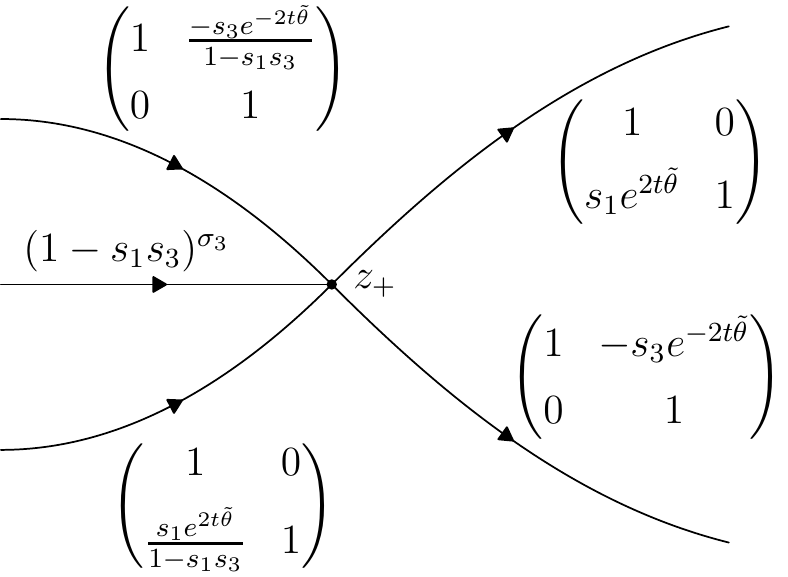}
\end{center}
\end{subfigure}
\caption{The contours $D(z_{\pm},\delta)\cap \Sigma_{T}$ for the RH problem satisfied by the local parametrix around the stationary points $z_{\pm}=\pm \frac{1}{2}$.}
\label{fig12}
\end{figure}
We define the functions $T^{(r)}(z)$ on the set $D(z_{+},\delta)\setminus \Sigma_{T}$, by 
\begin{align}\label{t-r-def} 
\hspace{-10pt}T^{(r)}(z)\!:=\!\beta(z)^{\sigma_{3}}\!\left(\!\frac{-h_{1}}{s_{3}}\!\right)^{\!\!-\frac{\sigma_{3}}{2}}\!\!\!\!e^{\frac{it\sigma_{3}}{3}}2^{-\frac{\sigma_{3}}{2}}\!\!
\begin{pmatrix} \sqrt{t}\zeta(z) & \!\!1 \\[3pt] 1& \!\!0\end{pmatrix}\!Z(\sqrt{t}\zeta(z)) e^{t\tilde\theta(z)\sigma_{3}}\!
\left(\!\frac{-h_{1}}{s_{3}}\!\right)^{\!\!\frac{\sigma_{3}}{2}}
\end{align}
where $\zeta(z)$ is a biholomorphic map given by the formula \eqref{d-zeta} and 
\begin{align*}
\beta(z) := \left(\sqrt{t}\zeta(z)\frac{z+1/2}{z-1/2}\right)^{\nu},
\end{align*}
with the branch cut chosen such that $\mathrm{arg}\,z\in(-\pi/2,\pi/2)$. 
\begin{proposition}\label{l-par}
The function $T^{(r)}(z)$ is a solution of the following RH problem.\\[3pt]
\noindent\makebox[5.5mm][l]{$(a)$}\parbox[t][][t]{121mm}{The function $T^{(r)}(z)$ is analytic in $D(z_{+},\delta)\setminus\Sigma_T$.}\\[3pt]
\noindent\makebox[5.5mm][l]{$(b)$}\parbox[t][][t]{121mm}{On the contour $\Sigma^{+}_{T} := D(z_{+},\delta)\cap \Sigma_{T}$ the function $T^{(r)}(z)$ satisfies the same jump conditions as $T(z)$ (see right diagram of Figure \ref{fig12}).}\\[3pt]
\noindent\makebox[5.5mm][l]{$(c)$}\parbox[t][][t]{121mm}{As $z\to z_{+}$, the function $T^{(r)}(z)$ is bounded.}\\[3pt]
\noindent\makebox[5.5mm][l]{$(d)$}\parbox[t][][t]{121mm}{The following asymptotic relation is satisfied
\begin{align}\label{asympt-1}
T^{(r)}(z)N(z)^{-1} =\! \begin{pmatrix} 1+\frac{\nu(\nu+1)}{2t\zeta(z)^{2}} & \frac{-\nu s_{3}}{h_{1}}e^{\frac{2it}{3}}\!\frac{\beta(z)^{2}}{t^{1/2}\zeta(z)} \\[3pt] 
\frac{-h_{1}}{s_{3}}e^{-\frac{2it}{3}}\!\frac{\beta(z)^{-2}}{t^{1/2}\zeta(z)} & 1-\frac{\nu(\nu-1)}{2t\zeta(z)^{2}} \end{pmatrix} + O(t^{-\frac{3}{2}}), \quad t\to\infty,
\end{align}
uniformly for $z\in \partial D(z_{+},\delta)$. }
\end{proposition}
The proof of the above proposition can be actually found in \cite[Section 9.4]{MR2264522}, where there was proved that the function $T^{(r)}(z)$ satisfies conditions $(a)-(c)$ and furthermore the following weaker asymptotic relation holds
\begin{align}\label{asym-11-22}
T^{(r)}(z)N(z)^{-1} = \begin{pmatrix} 1 & -\frac{\nu s_{3}}{h_{1}}e^{\frac{2it}{3}}\frac{\beta^{2}(z)}{\sqrt{t}\zeta(z)} \\[3pt]
-\frac{h_{1}}{s_{3}}e^{-\frac{2it}{3}}\frac{\beta^{-2}(z)}{\sqrt{t}\zeta(z)} & 1\end{pmatrix} + O(t^{-1}),\quad t\to\infty,
\end{align}
uniformly for $z\in\partial D(z_{+},\delta)$. Although the relation \eqref{asym-11-22} is insufficient to prove the asymptotic \eqref{asym-1a}, it can be readily improved to \eqref{asympt-1}, by an application of more accurate asymptotic expansions of the parabolic cylinder functions. To make our paper self-contained we shortly perform these calculations and consequently we determine the subsequent coefficient located by $t^{-1}$. \\[5pt]
\noindent{\em Proof of Proposition \ref{l-par}.} The results of \cite[Section 9.4]{MR2264522} (see also \cite[Section 3.4]{MR3670014}) say that the function $T^{(r)}(z)$ satisfies conditions $(a)-(c)$ and hence, it remains to derive the relation \eqref{asympt-1}. To this end, we consider the function $\xi(z):=\sqrt{t}\zeta(z)$, which clearly satisfies the asymptotic
\begin{align}\label{xi-asym}
|\xi(z)| = |\sqrt{t}\zeta(z)|= O(t^{1/2}),\quad t\to\infty, 
\end{align}
uniformly for $z\in\partial D(z_{+},\delta)$. Using the asymptotic expansions of the parabolic cylinder functions (see e.g. \cite{ba-er}, \cite{MR0232968}), we infer that
\begin{align*}
Z(z)=\frac{z^{-\frac{\sigma_{3}}{2}}}{\sqrt{2}}\left(\!\begin{pmatrix} 1 & 1 \\[3pt] 1 & -1\end{pmatrix} + \begin{pmatrix} \frac{(\nu+1)(\nu+2)}{2z^{2}} & -\frac{\nu(\nu-1)}{2z^{2}} \\[3pt] \frac{(\nu+1)(\nu-2)}{2z^{2}} & \frac{\nu(\nu+3)}{2z^{2}} \end{pmatrix} + \tilde Z(z)\right)e^{(\frac{1}{4}z^{2}-(\nu+\frac{1}{2})\ln z)\sigma_{3}},
\end{align*}
where the remainder term satisfies $\|\tilde Z(z)\|= O(z^{-4})$ as $z\to\infty$. This enables us to write the following decomposition 
\begin{align}\label{eq-p-n}
\begin{pmatrix} \sqrt{t}\zeta(z) & 1 \\[3pt] 1& 0\end{pmatrix}Z(\sqrt{t}\zeta(z)) = 2^{-1/2}(K_{1}(z) + K_{2}(z) + \tilde K(z))e^{(\frac{1}{4}z^{2}-\nu \ln z)\sigma_{3}},
\end{align}
where the component terms have the form 
\begin{align*}
K_{1}(z) := \begin{pmatrix} 2& 0 \\[3pt] \xi(z)^{-1}& 1\end{pmatrix}, \quad 
K_{2}(z) := \begin{pmatrix} \frac{\nu(\nu+1)}{\xi(z)^{2}} & \frac{2\nu}{\xi(z)} \\[3pt] \frac{(\nu+1)(\nu+2)}{2\xi(z)^{3}} & -\frac{\nu(\nu-1)}{2\xi(z)^{2}} \end{pmatrix}
\end{align*}
and furthermore
\begin{equation*}
\begin{gathered}
\tilde K_{11}(z) := \tilde Z_{11}(\xi(z)) + \tilde Z_{21}(\xi(z)), \quad \tilde K_{12}(z) :=  \xi(z)(\tilde Z_{12}(\xi(z))+\tilde Z_{22}(\xi(z))), \\
\tilde K_{21}(z) := \tilde Z_{11}(\xi(z))\xi(z)^{-1},\quad \tilde K_{22}(z) := \tilde Z_{12}(\xi(z)).
\end{gathered}
\end{equation*}
In view of the asymptotic \eqref{xi-asym}, we can write
\begin{equation}\label{equa-i-3}
K_{1}(z) + K_{2}(z)+ \tilde K(z)= L_{1}(z) + \tilde L(z)
\end{equation}
where we assume that
\begin{align*}
L_{1}(z):=\begin{pmatrix} 2+\frac{\nu(\nu+1)}{\xi(z)^{2}} & \frac{2\nu}{\xi(z)} \\[5pt] \frac{1}{\xi(z)} & 1-\frac{\nu(\nu-1)}{2\xi(z)^{2}} \end{pmatrix}
\end{align*}
and the remainder term $\tilde L(z)$ have the following asymptotic behavior as $t\to\infty$:
\begin{equation*}
\begin{aligned}
\tilde L_{11}(z) &= \tilde Z_{11}(\xi(z)) + \tilde Z_{21}(\xi(z)) = O(t^{-2}), \\ 
\tilde L_{12}(z) &=  \xi(z)(\tilde Z_{12}(\xi(z))+\tilde Z_{22}(\xi(z))) = O(t^{-3/2}), \\
\tilde L_{21}(z) &= \tilde Z_{11}(\xi(z))\xi(z)^{-1} + (\nu+1)(\nu+2)\xi(z)^{-3}/2 = O(t^{-3/2}), \\
\tilde L_{22}(z) &= \tilde Z_{12}(\xi(z)) = O(t^{-2}), 
\end{aligned}
\end{equation*}
In particular, we have $\|\tilde L(z)\| = O(t^{-\frac{3}{2}})$ as $t\to\infty$, uniformly for $z\in\partial D(z_{+},\delta)$. By \eqref{stokes-11bb} and \eqref{stokes-11bbc}, we know that $\nu$ is the real number and consequently
\begin{align}\label{eq-11-bb}
\beta^{\pm 1}(z) = O(1),\quad t\to\infty
\end{align}
uniformly for $z\in \partial D(z_{+},\delta)$. Combining this with the following equality
\begin{align}\label{eq-11-aa}
\!e^{(\frac{1}{4}\xi(z)^{2}-\nu\ln\xi(z))\sigma_{3}}e^{t\tilde\theta(z)\sigma_{3}}N(z)^{-1} \!= e^{-\frac{it\sigma_{3}}{3}}\xi(z)^{-\nu\sigma_{3}}N(z)^{-1} \!= e^{-\frac{it\sigma_{3}}{3}}\beta(z)^{-\sigma_{3}}
\end{align}
we obtain the following asymptotic relation as $t\to\infty$:
\begin{equation}
\hspace{-5pt}\begin{aligned}\label{equa-i-2}
&\frac{\beta(z)^{\sigma_{3}}}{\sqrt{2}}\left(\!\frac{-h_{1}}{s_{3}}\!\right)^{\!\!-\frac{\sigma_{3}}{2}}e^{\frac{it\sigma_{3}}{3}}2^{-\frac{\sigma_{3}}{2}}\tilde L(z)e^{(\frac{1}{4}\xi(z)^{2}-\nu\ln\xi(z))\sigma_{3}}e^{t\tilde\theta(z)\sigma_{3}}N(z)^{-1}\left(\!\frac{-h_{1}}{s_{3}}\!\right)^{\!\!\frac{\sigma_{3}}{2}} \\
& \quad = \frac{\beta(z)^{\sigma_{3}}}{\sqrt{2}}\left(\!\frac{-h_{1}}{s_{3}}\!\right)^{\!\!-\frac{\sigma_{3}}{2}}e^{\frac{it\sigma_{3}}{3}}2^{-\frac{\sigma_{3}}{2}}\tilde L(z)e^{-it\sigma_{3}/3}\left(\!\frac{-h_{1}}{s_{3}}\!\right)^{\!\!\frac{\sigma_{3}}{2}}\beta(z)^{-\sigma_{3}} = O(t^{-\frac{3}{2}})
\end{aligned}
\end{equation}
uniformly for $z\in\partial D(z_{+},\delta)$. Then, using \eqref{eq-11-bb} and \eqref{eq-11-aa} once again, we have
\begin{equation*}
\begin{gathered}
\frac{\beta(z)^{\sigma_{3}}}{\sqrt{2}}\left(\frac{-h_{1}}{s_{3}}\right)^{-\frac{\sigma_{3}}{2}}e^{\frac{it\sigma_{3}}{3}}2^{-\frac{\sigma_{3}}{2}}L_{1}(z)e^{(\frac{1}{4}\xi(z)^{2}-\nu\ln\xi(z))\sigma_{3}}e^{t\tilde\theta(z)\sigma_{3}}N(z)^{-1}\left(\frac{-h_{1}}{s_{3}}\right)^{\frac{\sigma_{3}}{2}}\\[5pt]
= \begin{pmatrix} 1+\frac{\nu(\nu+1)}{2\xi(z)^{2}} & \frac{-\nu s_{3}}{h_{1}}e^{\frac{2it}{3}}\frac{\beta(z)^{2}}{\xi(z)} \\[7pt] 
\frac{-h_{1}}{s_{3}}e^{-\frac{2it}{3}}\frac{\beta(z)^{-2}}{\xi(z)} & 1-\frac{\nu(\nu-1)}{2\xi(z)^{2}} \end{pmatrix}
\end{gathered}
\end{equation*}
which together with \eqref{eq-p-n}, \eqref{equa-i-3} and \eqref{equa-i-2}, gives the desired asymptotic \eqref{asympt-1}
and the proof of the proposition is completed. \hfill $\square$ \\

Using the symmetry of the contour $\Sigma_{T}$ we can define the local parametrix $T^{(l)}(z)$ around the stationary point $z_{-}=-1/2$ as
\begin{align}\label{t-l-def}
T^{(l)}(z):=\sigma_{2}T^{(r)}(-z)\sigma_{2},\quad z\in D(z_{-},\delta)\setminus \Sigma_{T}.
\end{align}
In the following proposition we provide analogous improvement of the asymptotic behavior of $T^{(l)}(z)$ as $z\to\infty$.
\begin{proposition}\label{l-par2}
The function $T^{(l)}(z)$ is a solution of the following RH problem.\\[3pt]
\noindent\makebox[5.5mm][l]{$(a)$}\parbox[t][][t]{121mm}{The function $T^{(l)}(z)$ is analytic in $D(z_{-},\delta)\setminus\Sigma_T$.}\\[3pt]
\noindent\makebox[5.5mm][l]{$(b)$}\parbox[t][][t]{121mm}{On the contour $\Sigma^{-}_{T} = D(z_{-},\delta)\cap \Sigma_{T}$ the function $T^{(l)}(z)$ satisfies the same jump conditions as $T(z)$ (see left diagram of Figure \ref{fig12}).}\\[3pt]
\noindent\makebox[5.5mm][l]{$(c)$}\parbox[t][][t]{121mm}{As $z\to z_{-}$, the function $T^{(l)}(z)$ is bounded.}\\[3pt]
\noindent\makebox[5.5mm][l]{$(d)$}\parbox[t][][t]{121mm}{The following asymptotic relation is satisfied
\begin{align}\label{asympt-2}
\!\!T^{(l)}(z)N(z)^{-1} = \begin{pmatrix}  1-\frac{\nu(\nu-1)}{2t\zeta(-z)^{2}} & \frac{h_{1}}{s_{3}}e^{-\frac{2it}{3}}\frac{\beta(-z)^{-2}}{t^{1/2}\zeta(-z)} \\[7pt] \frac{\nu s_{3}}{h_{1}}e^{\frac{2it}{3}}\frac{\beta(-z)^{2}}{t^{1/2}\zeta(-z)} & 1+\frac{\nu(\nu+1)}{2t\zeta(-z)^{2}} \end{pmatrix} + O(t^{-\frac{3}{2}}), \ \ t\to\infty
\end{align}
uniformly for $z\in \partial D(z_{-},\delta)$. }
\end{proposition}
\proof From the definition of the function $T^{(l)}(z)$ and Proposition \ref{l-par2} it follows that $T^{(l)}(z)$ satisfies the points $(a)-(c)$. To prove that the point $(d)$ holds true, let us observe that $\sigma_{2}N(z)^{-1} = N(-z)^{-1}\sigma_{2}$ and therefore
\begin{align}\label{equa-pp-kk}
T^{(l)}(z)N(z)^{-1}=\sigma_{2}T^{(r)}(-z)\sigma_{2}N(z)^{-1} = \sigma_{2}T^{(r)}(-z)N(-z)^{-1}\sigma_{2}. 
\end{align}
Combining \eqref{equa-pp-kk} and \eqref{asympt-1} we have the following asymptotics as $t\to\infty$:
\begin{align*}
T^{(l)}(z)N(z)^{-1} & = \sigma_{2}\begin{pmatrix} 1+\frac{\nu(\nu+1)}{2\xi(-z)^{2}} & \frac{-\nu s_{3}}{h_{1}}e^{\frac{2it}{3}}\frac{\beta(-z)^{2}}{\xi(-z)} \\[7pt] \frac{-h_{1}}{s_{3}}e^{-\frac{2it}{3}}\frac{\beta(-z)^{-2}}{\xi(-z)} & 1-\frac{\nu(\nu-1)}{2\xi(-z)^{2}} \end{pmatrix}\sigma_{2} + O(t^{-\frac{3}{2}})\\[5pt]
&  = \begin{pmatrix}  1-\frac{\nu(\nu-1)}{2\xi(-z)^{2}} & \frac{h_{1}}{s_{3}}e^{-\frac{2it}{3}}\frac{\beta(-z)^{-2}}{\xi(-z)} \\[7pt] \frac{\nu s_{3}}{h_{1}}e^{\frac{2it}{3}}\frac{\beta(-z)^{2}}{\xi(-z)} & 1+\frac{\nu(\nu+1)}{2\xi(-z)^{2}} \end{pmatrix} + O(t^{-\frac{3}{2}}),
\end{align*}
which establishes \eqref{asympt-2} and the proof of the proposition is completed. \hfill $\square$ \\

\noindent {\bf Acknowledgements.} 
%The authors would like to thank the referees for their helpful comments and remarks. 
The second author is supported by the MNiSW Iuventus Plus Grant no. 0338/IP3/2016/74.

\parindent = 0 pt

\end{document}